\documentclass[a4paper, 12pt, oneside]{article}
\usepackage{graphicx}
\usepackage{amsmath}
\usepackage{amssymb}
\usepackage{amsfonts}
\usepackage{amsthm}
\usepackage{longtable}
\usepackage{multirow}
\usepackage{color}

\newtheorem{lemma}{Lemma}[section]
\newtheorem{theorem}[lemma]{Theorem}
\newtheorem{proposition}[lemma]{Proposition}
\newtheorem{corollary}[lemma]{Corollary}

\renewcommand{\Gamma}{\varGamma}
\renewcommand{\epsilon}{\varepsilon}
\renewcommand{\bar}{\overline}

\newcommand{\type}{\mathcal{T}}

\begin{document}
\title{Vertex-transitive maps with Schl\"afli type $\{3, 7\}$}
\author{Daniel Pellicer}

\maketitle

\begin{abstract}
Among all equivelar vertex-transitive maps on a given closed surface $S$,
the automorphism groups of maps with Schl\"afli types $\{3, 7\}$ and
$\{7, 3\}$ allow the highest possible order. We describe a procedure to transform all such
maps into $1$- or $2$-orbit maps, whose symmetry type has been previously studied. In so doing we provide
a procedure to determine all vertex-transitive maps with Schl\"afli type $\{3, 7\}$ which
are neither regular or chiral. We determine all such maps on surfaces with Euler characteristic $-1 \ge \chi \ge -40$.
\end{abstract}


\section{Introduccion}

   Highly symmetric maps on compact surfaces have attracted attention in recent years. Several
efforts have been made to classify such maps in surfaces with given Euler characteristics.

   Regular maps (maps with maximal symmetry by reflections) and chiral maps
(maps with maximal symmetry by rotations, but no symmetry by
reflections) have been the most studied.
   Regular and chiral maps with Euler
characteristic $-1 \ge \chi \ge -200$ have been classified by Conder in \cite{conderatlas}. If
we additionally require the so-called diamond condition, regular maps whose
automorphism group contains at most 2000 elements were classified
as rank $3$ abstract polytopes by Hartley in \cite{sma}, with the exception of those whose
automorphism group contains $1024$ or $1536$ elements.

   Edge-transitive maps (maps whose automorphism group acts transitively on their edges) were
studied and classified into $14$ symmetry types by several authors
in \cite{GW97} and \cite{STW2001}. An atlas of each of the types of edge-transitive maps up
to certain Euler characteristic depending on the type was presented by Orbani\'c in
\cite{alenDoktorat}.

   Vertex-transitive maps have a wider variety of symmetry types than
edge-transitive maps. This makes considerably more
complicated to determine all such maps on a given surface. A combinatorial approach to
attack the problem of a classification of vertex-transitive maps on surfaces with a given
Euler characteristic is outlined in \cite[Section 10]{Pellicer}. Another approach to this classification, is given by Karab\'a\v{s} and Nedela in \cite{karabas} and \cite{karabas2},
where extra properties are required to the maps. This approach involves actions of groups on
surfaces. They also present in \cite{karabas} an atlas of the so-called Archimedean
solids on an orientable surface with Euler characteristic $-2$.

   We may restrict ourselves to determine all vertex-transitive maps satisfying the extra
condition that all its faces have the same co-degree $p$. The automorphism groups of these
maps on a surface with Euler characteristic $-m$ contain at most $4mpq/(pq - 2p -2q)$ elements,
where $q$ is the degree of the vertices. This bound is achieved if and only if
the map is regular (see \cite[Section 1]{conderdob}). An easy calculation shows that
the maps that maximize this bound are those where $p=3$ and $q=7$, or $p=7$ and $q=3$.
Furthermore, if we eliminate these cases, the bound is reduced by a factor of
$12/21$.

   The most common method to produce an atlas of maps with symmetries consists in an exhaustive
search of possible automorphism groups with distinguished generators from which the map can be
recovered. In this context, the reduction in the bound described above may play a significant role
in classifying vertex-transitive maps with faces of the same size.

   As a consequence of the notions of maps on surfaces and automorphisms, there is a variety of
algebraic, geometric, topologic and combinatorial techniques to work with these concepts.
In this paper we combine the technique mentioned in the paragraph above with operations on maps,
which are consistent transformations on all fundamental regions of the map under a certain automorphism group. One advantage
of using operations is that we naturally obtain local properties of the map and relations to other maps, that may be difficult to spot otherwise. Additionally,
these operations can be applied to higher dimentional objets, when the flags are no longer triangles, but $n$-simplexes.

   In Section \ref{sec:def} we recall some definitions and basic results about symmetric maps.
We define operations on maps in Section \ref{sec:oper}. These operations are used
in Section \ref{section:3773} to transform vertex-transitive triangulations with seven triangles
around each vertex into maps whose symmetry type has been previously studied. This reduction
process can be reversed to construct all such maps from maps with a simple symmetry type.
We also show that vertex-transitive maps all whose faces are heptagons, and all whose vertices
are trivalent, must necessarily be either regular or chiral. In Section \ref{sec:exam}
we construct three vertex-transitive maps to illustrate the procedure described in Section
\ref{section:3773}. Finally, in an appendix we list all vertex-transitive triangulations with seven
triangles around each vertex in surfaces with Euler characteristic $-1 \ge \chi \ge -40$.

\section{Definitions}\label{sec:def}

   Throughout, a {\em map} is a 2-cell embedding of a finite multigraph (we allow multiple edges) $G$ on a closed
surface $S$ without boundary. The connected components of $S$ after removing $G$
are called {\em faces} and are homeomorphic to discs. The reader is referred to
\cite[Section 3]{Pellicer} for details. Equivalent definitions using different approaches are given for example in \cite{conderatlas}, \cite{JonesSin}, \cite{alenasiadan}, \cite{Skoviera} and \cite{vince}.

   A map $M$ on a surface $S$ naturally induces a triangulation on $S$ where the
vertices of the triangles are the vertices, midpoints of edges and centres of faces
of $M$, with all triangles containing one vertex of each type. These triangles are
henceforth called {\em flags}. We say that two flags are $0$-adjacent (resp.
$1$- and $2$-adjacent) whenever they share a line segment between the midpoint of the edge
and the centre of the face (resp. a line segment between the vertex and the centre of the
face, and a line segment between the vertex and the midpoint of the edge).

   We define the {\em flag graph} of a map $M$ on a surface $S$
as the graph whose vertices are the
flags of $M$ and two vertices are joined by an edge labelled $i$ whenever the
corresponding flags are $i$-adjacent ($i = 0, 1, 2$). Note that the flag graph of
$M$ allows a $2$-cell embedding on $S$. It is easy to see that
the flag graph of a map is always
a connected graph without loops where the edges of each color form a matching. %In
%general the flag graph of a map may contain multiple edges, however this will not be the
%case in this paper.

   We say that a map is {\em equivelar} whenever the {\em co-degree} of all faces
(number of edges around the faces) is a fixed number $p$ and the degree of all
vertices is a fixed number $q$.

   An {\em automorphism} of a map $M$ is an automorphism of its graph that can be
extended to an homeomorphism of the surface. We shall denote the group consisting
of all automorphisms of $M$ by $\Gamma(M)$.

   We now introduce the type graph of a map, which
is a variation of the Delaney-Dress graph described in
\cite{Dr1} and \cite{Dr2}, with semi-edges instead of loops.

   The {\em type graph} $\type(M)$ of a map $M$ is the semi-graph (we may allow
semi-edges) with edges and semi-edges labelled in $\{0, 1, 2\}$ defined
as follows. The vertices of $\type(M)$ are the flag-orbits of $M$ under $\Gamma(M)$,
with two of them adjacent by an edge labelled $i$ whenever the flags in one of the orbits
are $i$-adjacent to the flags in the other. Whenever a flag in an orbit
$\mathcal{O}$ is $i$-adjacent to another flag in $\mathcal{O}$ we attach a semi-edge
labelled $i$ to the vertex corresponding to $\mathcal{O}$. Note that the flag graph
and the type graph of a map are $3$-valent with one edge (or semi-edge)
of each color incident to each vertex. In this sense, the
type graph of $M$ is the quotient of the flag graph by $\Gamma(M)$. We shall refer
to all maps with type graph $\type$ as {\em maps of type} $\type$.

   As opposite to the flag graphs, we shall consider type graphs with multiple edges.

   Whenever $\type$ is a bipartite type graph (not containing semi-edges or odd cycles),
every map $M$ of type $\type$ must lie on an orientable surface. In fact,
if $(f_1, \dots, f_s=f_1)$ is a sequence
of flags such that $f_i$ is adjacent to $f_{i+1}$ for $i = 1, \dots, s-1$ then
$s$ is odd, and therefore, the sequence involves an even number of changes
of flags. If $\type$ is not bipartite then $M$ may lie either
on an orientable or on a non-orientable surface.

   In what follows, given the type graph $G$ of a map we denote by $G_i$
the subgraph of $G$ obtained from $G$ by deleting all edges labelled $i$
($i=0, 1, 2$).

\begin{figure}
\begin{center}
\includegraphics[width=8cm, height=2.5cm]{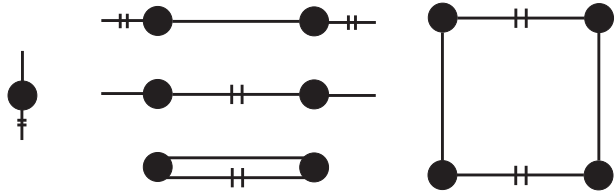}
\caption{Connected components of $G_1$}\label{fig:graph02}
\end{center}
\end{figure}

   We note that if $G$ is the type graph of a map, then each connected
component of $G_1$ is isomorphic to one of the ones showed in Figure
\ref{fig:graph02}.

   The number of orbits of vertices (resp. edges, faces) of a map with type graph $G$
coincides with the number of connected components of $G_0$ (resp. $G_1$, $G_2$).
This implies that {\em edge-transitive maps} (maps where the automorphism group
acts transitive on its edges) can have at most $4$-orbits on flags, however there is
not such a bound for the number of orbits of flags for vertex-transitive maps
(maps where the automorphism group acts transitive on its vertices).

   A map $M$ whose automorphism group induces $k$ flag orbits is called a
{\em $k$-orbit map}. In this case $|\Gamma(M)| = |\mathcal{F}l(M)|/k$, where
$\mathcal{F}$ denotes the set of flags of $M$.

   We shall follow the notation of \cite{alenasiadan} and
denote $1$-orbit maps as {\em regular} maps. A {\em chiral} map is a
$2$-orbit map such that adjacent flags belong to distinct orbits. Note that in \cite{conderdob}
regular maps are referred to as {\em reflexible}.

   The seven types of $2$-orbit abstract polyhedra are described in \cite{isabelspaper} and we
shall largely follow its notation. In fact, every finite abstract polyhedron can be seen as a map on a surface, and all maps on surfaces satisfying the so-called diamond condition are also abstract polyhedra. The theory developed in \cite{isabelspaper} is valid also for maps which are not polyhedra if we ignore the diamond and intersection conditions.

\begin{figure}
\begin{center}
\includegraphics[width=7cm, height=2cm]{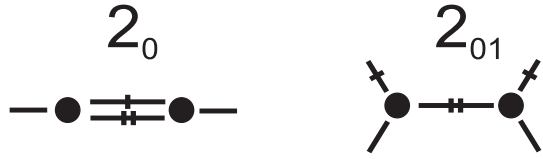}
\caption{Types graphs $2_0$ and $2_{01}$}\label{fig:2orbits}
\end{center}
\end{figure}

    Among the seven types of $2$-orbit maps, of particular interest for us are the types denoted by
$2_0$ and $2_{01}$ (see Figure \ref{fig:2orbits} for their type graphs). These types are respectively called $2^*ex$ and $2^*$ in \cite{GW97} and \cite{STW2001}, where they consider only
maps on orientable surfaces. Every flag of a map of type $2_0$ is $0$-adjacent to a flag
in the same orbit, and is $1$- and $2$-adjacent to flags in the other orbit.
Similarly, every flag of a map of type $2_{01}$ is $0$- and $1$-adjacent to
flags in the same orbit, and is $2$-adjacent to a flag in the other orbit.

   The automorphism group of any map of type $2_0$ is generated by an involution
$\rho_0$ and a non-involution $\sigma_{12}$, respectively denoted by $\lambda_e$
and $\sigma_f$ in \cite{STW2001}.
   The automorphism group of any map of type $2_{01}$ is
generated by three involutions $\rho_0$, $\rho_1$ and $\rho_{212}$ (the latter denoted by $\alpha_{212}$ in \cite{isabelspaper}), respectively
denoted by $\lambda_e$, $\theta_{uf}$ and $\theta_{uf} \sigma_u^{-2}$
in \cite{STW2001}.

   Given a type $\type$ of $k$-orbit maps with corresponding type graph $G$
we say that a map $M$ is {\em $\type$-admissible} if it admits a labelling of
its flags with labels in the vertex set of $G$ such that
\begin{itemize}
   \item flags of $M$ labelled $x$ are $i$-adjacent to flags of $M$ labelled
$y$ if and only if there is an edge labelled $i$ between the vertices $x$ and
$y$ of $G$.
   \item all flags of $M$ with the same label belong to the same orbit.
\end{itemize}
Note that $\type$-admissible maps correspond to $\type$-regular maps in
\cite{breda}.

   Clearly, if $\type$ is a type of $k$-orbit maps then every $\type$-admissible
map is a $k'$-orbit map, for some divisor $k'$ of $k$. In particular, if $\type$ is
a type of $2$-orbit maps then every $\type$-admissible map is either of type
$\type$ or regular.

\section{Operations}\label{sec:oper}

   In Sections \ref{section:3773} and \ref{sec:exam} we shall require some operations
on maps described next.

   An {\em operation} $Oper$ applied to the map $M$ on the surface $S$ yields
a map $Oper(M)$ obtained from $M$ by the following steps.

\begin{itemize}
   \item[(a)] Divide $S$ in a set of regions $R_1, \dots, R_k$ isomorphic to discs.
      Each region inherits a portion of a graph from the map $M$.
   \item[(b)] Modify the portion of graph on each $R_i$.
   \item[(c)] Identify the boundaries of the regions $R_i$ in some way determined by
      $Oper$ to form a surface $S'$. $Oper(M)$ will be the induced map on $S'$.
\end{itemize}

   Next we consider the operations dual, Petrial and truncation, as well as
two specific operations, called collapsing and rebelting, that can only be applied to
certain families of maps. In all these operations each region consists of a
collection of flags.

\subsection{Dual}

   The {\em dual} operation consists of reinterpreting the centres of faces as vertices
and the vertices as centres of faces preserving the same number of edges. This can
be done by considering the same set of flags, but interchanging $0$- and
$2$-adjacencies. Alternatively, we may define it in the following way.

   Let $M$ be a map on a surface $S$. For every face $f$ of $M$ consider the
region $R_f$ determined by all flags contained in $f$. On each region $R_f$, the
portion of the graph consists of a $p$-cycle around its boundary, where $p$ is
the co-degree of the face. Replace this portion of graph with a vertex in the
interior of the face with $p$ semi-edges which join the midpoint of all
original edges to the centre of $f$. Finally, identify the boundaries of the
regions $R_f$ in the way they were identified before to obtain the dual
$Du(M)$ of $M$. Each edge of $Du(M)$ intersects precisely one edge of $M$ and
therefore the edges of these two maps are in a natural one-to-one correspondence.

   Note that $Du(M)$ lies on the same surface of $M$. Furthermore, if $M$ is
equivelar with Schl\"afli type $\{p, q\}$ then $Du(M)$ is equivelar with
Schl\"afli type $\{q, p\}$.

   It is easy to see that $\Gamma(M) \cong \Gamma(Du(M))$ and, therefore, $M$ is a
$k$-orbit map if and only if $Du(M)$ is a $k$-orbit map. Clearly $Du(Du(M)) = M$
for every map $M$.

\subsection{Petrial}

   We define the {\em Petrie polygons} (also called {\em Petrie paths})
of a map $M$ as the faces of $Pe(M)$. The Petrie polygons of $M$
can be visualized in $M$ as zigzags with the property that every two
consecutive edges of a Petrie polygon are two consecutive edges of a face; but
three consecutive edges of a Petrie polygon belong to the same face only when
the vertex between two of the edges has degree $2$.

   Often the Petrial of a map $M$ is defined in the map with the same vertex and edge set as $M$
with faces given the Petrie polygons of $M$. Alternatively we can define it in the following way.

   For each vertex $v$, consider the region $R_v$ consisting of all
flags around $v$. Each line-segment between two centres
$c_1$ and $c_2$ of faces sharing an edge $e = \{v_1, v_2\}$ in $M$ is shared by the
two regions $R_{v_1}$ and $R_{v_2}$. For each such line-segment, identify the corresponding
two regions reversing the original local orientation, that is, identifying the
line segment in such a way that the point $c_1$ on the region $R_{v_1}$ is
identified with the point $c_2$ of the region $R_{v_2}$,
and the point $c_2$ on the region $R_{v_1}$ with the point $c_1$ of the region
$R_{v_2}$. The resulting map is the {\em Petrial} of the map $M$ and is denoted by
$Pe(M)$. Note that $Pe(M)$ contains the same vertex and edge set as $M$.

   Since the underlying graph of $M$ is invariant
under the Petrie operation, $M$ has all vertices with degree $q$ if and only if
$Pe(M)$ also does. However, if $Pe(M)$ is equivelar it is not necessarily true that
$M$ is equivelar, but it is true that all Petrie polygons of $M$ have the same size.

   In general the Petrie operation changes the surface. In fact, the number $v$ of vertices
and the number of edges $e$ are preserved under this operation, but the number $f$ of faces
needs not be. Therefore the well-known surface invariant
{\em Euler characteristic} $v-e+f$ of a map $M$ may differ
with that of its Petrial. Furthermore, the underlying surface of $Pe(M)$ may
or may not be orientable independently to the orientability of the underlying surface of $M$.

   Clearly $Pe(Pe(M))=M$. Furthermore, if $M$ is a vertex-transitive
$k$-orbit map then $Pe(M)$ is also a vertex-transitive $k$-orbit map.

\subsection{Truncation}

   In \cite[Section 4.2]{alenasiadan} (see also \cite[Section 3]{Pellicer}) the
{\em truncation} operation is defined in three different ways. Here we present an
equivalent definition.

   For each vertex $v$ of $M$ consider the region $R_v$ consisting of all flags
around $v$. The portion of graph on each region consists of a vertex incident to
$q$ semi-edges, where $q$ is de degree of $v$. Replace $v$ by a $q$-cycle in such
a way that the new vertices lie in the interior of the old semi-edges and
identify the regions $R_v$ in the original way. The
resulting map is the truncation $Tr(M)$ of $M$. Note that the edges of $Tr(M)$ that
intersecting two different regions $R_v$ are in a one-to-one correspondence with
the edges of $M$. We call these edges {\em inherited edges} of $Tr(M)$.

   The truncation operation preserves the surface and yields exclusively 3-valent
maps. On the other hand, if $M$ is
a $k$-orbit map then $Tr(M)$ is either a $k$-, $(3k/2)$- or $(3k)$-orbit map
(see \cite[Proposition 4.5]{alenasiadan}). Note that $\Gamma(M)$ is a subgroup
of $Tr(M)$.

\subsection{Collapsing}

   Assume that there is an orbit $\mathcal{F}$ of triangular faces of
$M$ with the following two properties:
\begin{enumerate}
   \item every vertex is incident to either no faces in $\mathcal{F}$,
one face in $\mathcal{F}$, two faces in $\mathcal{F}$ sharing an edge,
or to three faces in $\mathcal{F}$ two of which share an edge,
   \item each triangle in $\mathcal{F}$ shares precisely one edge with
another face in $\mathcal{F}$,
   \item every triangle in $\mathcal{F}$ has three distinct vertices, and
any two triangles in $\mathcal{F}$ sharing an edge contain four
distinct vertices.
\end{enumerate}
The last item can be reworded as follows. The closure of the region
determined by any two triangles in $\mathcal{F}$ sharing an edge is
homeomorphic to a disc.

\begin{figure}
\begin{center}
\includegraphics[width=15cm, height=5cm]{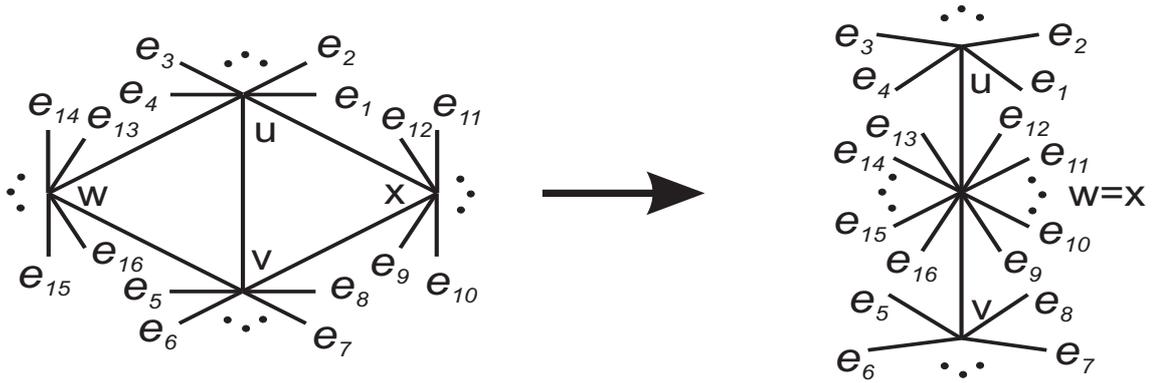}
\caption{Collapsing operation}\label{fig:detriangulating}
\end{center}
\end{figure}

   We define the {\em collapsing} operation with respect to $\mathcal{F}$
as follows. For each pair $F_1$ and
$F_2$ of faces in $\mathcal{F}$ sharing an edge $uv$, let $x$ be the vertex
in $F_1 \smallsetminus F_2$, and let $w$ be the vertex in $F_2 \smallsetminus
F_1$. Delete the triangles $F_1$ and $F_2$ by identifying the vertices $w$ and
$x$, the edges $ux$ and $uw$, and the edges $vx$ and $vw$
(see Figure \ref{fig:detriangulating}).

   In terms of the definition of operation above,
we may define the collapsing operation by considering the region $R_v$ for
each vertex $v$ as in the operations Petrial and truncation,
and eliminating from it the regions induced by all triangles
in $\mathcal{F}$. If the eliminated region corresponded to two adjacent triangles,
the semi-edges incident to $v$ bounding this region are identified. If the eliminated region
corresponded to only one triangle then the two semi-edges delimiting this region become
part of the boundary of the new region. The identification rule of the new regions is
inherited from that in $M$ except in the part of the boundary arising from the elimination
of a triangle $T_1$ in $\mathcal{F}$. By hypothesis $T_1$ is adjacent to precisely one
triangle $T_2$ in $\mathcal{F}$. Let $v'$ be the vertex in $T_2$ not contained in
$T_1$. Then the part of the boundary of $R_v$ corresponding to the semi-edges bounding
$T_1$ is identified to the part of the boundary of $R_{v'}$ arising from the
deletion of $T_2$.

   Note that Items $1$ and
$3$ guarantee that the collapsing operation preserves the surface.

\subsection{Rebelting}

\begin{figure}
\begin{center}
\includegraphics[width=8cm, height=1.5cm]{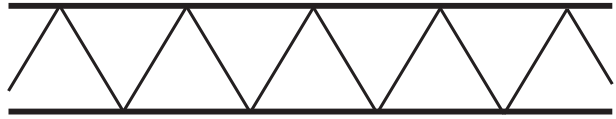}
\caption{Band of triangles}\label{fig:band}
\end{center}
\end{figure}

   The rebelting operation can only be applied to maps with vertices of degree $8$
whose faces are triangles with the properties that
\begin{itemize}
   \item some triangles are arranged in belts like the one in Figure \ref{fig:band}
such that no edge or face belongs to the intersection of two belts;
   \item the faces around every vertex are six triangles on belts and two triangles not belonging to belts, arranged in such a way that between the two triangles not belonging to belts there are precisely $3$ triangles on belts.
In other words, triangles which are not on belts share edges only with triangles on belts, and all triangles sharing edges with triangles on a belt are not on a belt.
\end{itemize}

\begin{figure}
\begin{center}
\includegraphics[width=8cm, height=4cm]{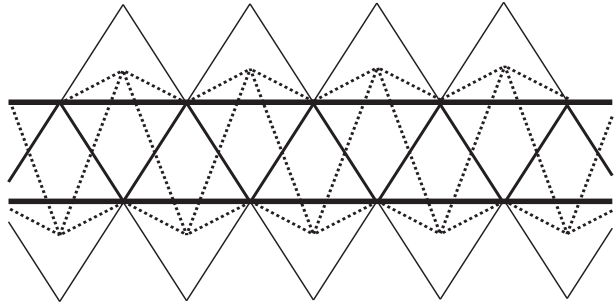}
\caption{Rebelting operation}\label{fig:rebelting}
\end{center}
\end{figure}

   The {\em rebelting operation} consists of substituting a belt of triangles by a
new belt with the same number of (bigger) triangles including one third of each triangle
adjacent to the original belt. In Figure \ref{fig:rebelting} the original (dark) belt
is substituted by the belt of dotted triangles. Note that the vertices of the original
triangulation become midpoints of edges, and the new vertices are the centres of triangles
not belonging to belts in the original triangulation.

   In terms of the definition of operation at the beginning of the section,
we define the collapsing operation as follows. For each triangle $t$ on a belt let $e_t$ be the edge of $t$ shared by $t$ and a triangle not on a belt.
We consider the region $R_t$ consisting of $t$ together with the two flags
around $e_t$ not belonging to $t$.
The portion of graph on each region consists of a triangle with precisely
two edges in the boundary of the region. We substitute this graph by a vertex $v_t$
with four semi-edges incident to it. The vertex $v_t$ is placed in the
point of the boundary corresponding to the centre of the triangle not belonging to a
belt in the original map. Two of the semi-edges join $v_t$ through the boundary of
$R_t$ with the two original vertices belonging to $t$ and to the triangle with
centre in $v_t$. The other ends of the remaining two semi-edges are the
midpoints of the two original edges in the boundary of $R_t$. We identify the regions
$R_t$ in the way they were identified in the original map. Note that the vertex $v_t$
appears in precisely three regions, and therefore, the resulting map has vertices
with degree $9$. Furthermore, all faces are triangles.

   Clearly the rebelting operation preserves the surface.

\section{$\{7, 3\}$ and $\{3, 7\}$ vertex-transitive maps}\label{section:3773}

   In this section we classify the vertex-transitive maps with Schl\"afli type
$\{3, 7\}$ which are neither regular nor chiral in three families. We also describe a
procedure to obtain them from regular or $2$-orbit maps by means of operations.
Additionally we prove that all vertex-transitive maps with Schl\"afli type $\{7, 3\}$
are either regular or chiral.

%   The following lemma is a direct consequence of the definition of Euler
%characteristic.

%\begin{lemma}\label{lemma:novertices}
%   The number of vertices of a map with Schl\"afli type
%$\{p, q\}$ on a surface with negative Euler characteristic $\chi$ is
%\[\frac{2p\chi}{2p-pq+2q}.\]
%\end{lemma}

%   In particular, a map with Schl\"afli type $\{3, 7\}$ on a surface of
%orientable genus $g \ge 2$ or nonorientable genus $g \ge 3$
%has $6 \chi$ vertices.

   We shall make use of the following lemmas.

\begin{lemma}\label{lemma:divides}
   Let $k$ be the number of orbits on flags of a vertex-transitive map with vertices of degree $q$. Then $k$ divides $2q$.
\end{lemma}

\begin{proof}
   The orbits of flags around any vertex $v$ induced by the stabilizer of $v$
have the same cardinality.
\end{proof}

\begin{lemma}\label{lemma:dividesequivelar}
   Let $M$ be an equivelar map with Schl\"afli type $\{p, q\}$.
Let $G$ be its type graph and let $C$ be a
connected component of $G_2$ (resp $G_0$). Then
\begin{itemize}
   \item if $C$ is a path, the number of vertices of $C$ divides $p$ (resp $q$),
   \item if $C$ is a cycle, the number of vertices of $C$ divides $2p$ (resp $2q$).
\end{itemize}
\end{lemma}

\begin{proof}
   The orbits represented in the vertices of $C$ correspond to the orbits
containing the flags in some face $F$ of $M$. Among the flags in
$F$, the stabilizer of $F$ induces orbits with the same cardinality. Therefore
the number of orbits of flags in $F$ must divide $2p$. Furthermore, if $C$ is a path,
then there is some orbit containing flags $i$-adjacent to flags in the same orbit for $i=0$ or $i=1$ and therefore,
the stabilizer of $F$ contains a reflection. Consequently,
each orbit $O$ of flags containing flags in $F$ contains
at least two flags contained on $F$, moreover, the
number of flags in $O \cap F$ is even. This implies the first item.

   A dual argument considering the flags around a vertex instead those in a face implies the lemma for $G_0$.
\end{proof}

   We are now ready to state our main result about maps with Schl\"afli type $\{p, 3\}$.

\begin{theorem}\label{prop:p7}
   All vertex-transitive maps with Schl\"afli type $\{p, 3\}$ with $p$ congruent
to $1$ or $5$ modulo $6$ are either regular or chiral.
\end{theorem}

\begin{proof}
   Let $p$ be congruent to $1$ or $5$ modulo $6$ and $G$ be the type graph of a vertex-transitive map $M$ with
Schl\"afli type $\{p, 3\}$.

   By Lemma \ref{lemma:divides}, the number of orbits on flags $M$ must be either
$1$, $2$, $3$ or $6$.
The dual version of Lemma \ref{lemma:dividesequivelar} implies that the number
of vertices on $G_0$
is either $1$ or $3$ if it is a path, or $2$ or $6$ if it is a cycle. Since $M$ is vertex-transitive, $G_0$ must be connected.
On the other hand, Lemma \ref{lemma:dividesequivelar}
implies that the number of vertices on each connected component of $G_2$
divides $2p$, that is, it is either $1$, $5$ or $2$. In
the latter case, it must be a cycle of length $2$.

   It is not hard to see that it is not possible to construct a graph
with $6$ vertices with the following three properties. The edges with labels $1$ and $2$
form a spanning cycle, the connected components induced by edges with labels $0$ and $2$
are among those in Figure \ref{fig:graph02}, and all connected components induced by
edges with labels $0$ and $1$ are either cycles of length $2$ or path with five vertices
(note that every vertex is adjacent to an edge labeled $1$ and an edge labeled $2$).

   Similarly it can be verified that it is not possible to construct a graph with
three vertices satisfying the following three conditions. The edges with labels
$1$ and $2$ form a spanning path, the connected components induced by edges with
labels $0$ and $2$ are among those in Figure \ref{fig:graph02}, and all
connected components induced by edges with labels $0$ and $1$ are either double edges
or isolated vertices.
Therefore any vertex-transitive
map with Schl\"afli type $\{p, 3\}$ has at most $2$ orbits on flags.

   Assume that $M$ has two orbits on flags. Then $G$
has only two vertices and an edge labelled $i$ between them for at least some $i$.
Lemma \ref{lemma:dividesequivelar}
implies that there must also be an edge with
label $i+1$ between the two vertices if $i < 2$, and an
edge with label $i-1$ between the two vertices if $i > 0$. This implies that
there are edges of all labels between the two vertices and hence the map is chiral.
Alternatively we could have argued that the remaining types of $2$-orbit maps
contain maps with either vertices with even valency, or faces with an even number
of edges (see \cite{isabelsthesis}).
\end{proof}

   The following corollary is a direct consequence of the previous theorem.

\begin{corollary}\label{prop:37}
   All vertex-transitive maps with Schl\"afli type $\{7, 3\}$ are either regular
or chiral.
\end{corollary}

   Now we turn our attention to vertex-transitive maps with Schl\"afli type
$\{3, 7\}$.

   Lemma \ref{lemma:divides} implies that any vertex-transitive equivelar map $M$
with Schl\"afli type $\{3, 7\}$ has either $1$, $2$, $7$ or $14$ orbits on flags.
Furthermore, similar arguments to those of the proof of Theorem
\ref{prop:p7} show that if $M$ has two orbits on flags it must be chiral. We first
eliminate the possibility of equivelar vertex-transitive maps with Schl\"afli type
$\{3, 7\}$ and $7$ orbits on flags.

\begin{proposition}\label{prop:7orbit}
   There are no $7$-orbit vertex-transitive equivelar maps with
Schl\"afli type $\{3, 7\}$.
\end{proposition}

\begin{proof}
   Assume to the contrary that such a map $M$ exists, then the type graph
$G$ of $M$ has $7$ vertices and $G_0$ is an alternating path with
$6$ edges. On the other hand, as a consequence of Lemma
\ref{lemma:dividesequivelar}, every connected component of $G_2$
is either a single vertex, two vertices joined by a double edge
(one edge of each color), an alternating
path with two edges, or an alternating cycle of length $6$.

   Let $v$ be the
vertex of $G$ incident to an edge labelled $1$ but to no edge labelled $2$.
Since the connected components induced by all edges with labels $0$ and $2$
are those in Figure \ref{fig:graph02}, and all vertices different form $v$ are
incident to an edge labelled $2$, $v$ cannot be incident to an edge labelled
$0$. Therefore the connected component of $G_2$ containing $v$
must be a path with two edges and
it must contain the vertex of $G$ not incident to an edge labelled $1$.
This fact, together with the allowed connected components of $G_1$ forces
another edge labelled $0$ forming an alternating square of edges labelled
$0$ and $2$, but this in turn induces a
connected component of $G_2$ different from the ones allowed.
\end{proof}

   We now describe all types of $14$-orbit maps containing vertex-transitive
maps with Schl\"afli type $\{3, 7\}$. Let $M$ be one such map and let
$G$ be its type graph.

\begin{figure}
\begin{center}
\includegraphics[width=12cm, height=4cm]{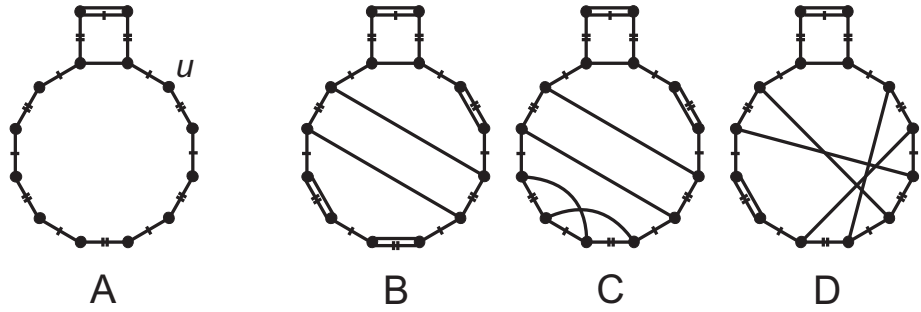}
\caption{Type graphs of $14$-orbit vertex-transitive maps with type
$\{3, 7\}$}\label{fig:14orbits}
\end{center}
\end{figure}

\begin{figure}
\begin{center}
\includegraphics[width=13.5cm, height=4cm]{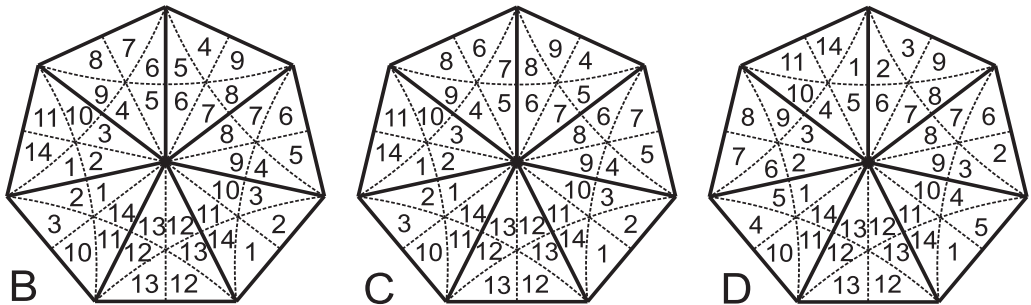}
\caption{Local configurations of $14$-orbit vertex-transitive maps with type
$\{3, 7\}$}\label{fig:14config}
\end{center}
\end{figure}

   Since $G_0$ is connected and has $14$ vertices, it must be an alternating
cycle. This implies that every vertex of $G$ is incident to an edge with label
$1$. Therefore every connected component of $G_2$ is either a double edge or
a $6$-cycle (see Lemma \ref{lemma:dividesequivelar}). By divisibility reasons
at least one connected component of $G_2$ must be a $2$-cycle and therefore
the graph in Figure \ref{fig:14orbits} (A) is a subgraph of $G$. An exhaustive
search for the choices of the vertex sharing an edge labeled $0$ with the
vertex labelled $u$ in this figure proves that the only possible type
graphs for vertex-transitive $14$-orbit maps with Schl\"afli type $\{3, 7\}$
are those in (B), (C) and (D) in Figure \ref{fig:14orbits}.
They correspond to the flag arrangements in Figure \ref{fig:14config}

\begin{figure}
\begin{center}
\includegraphics[width=13.5cm, height=4cm]{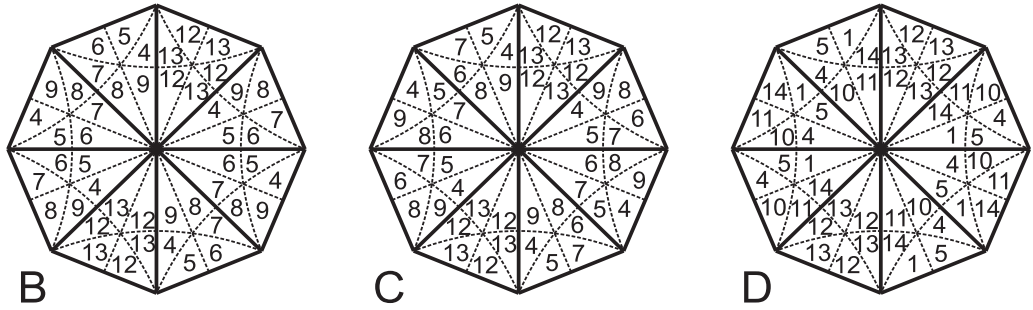}
\caption{Local configurations of $8$-orbit vertex-transitive maps with type
$\{3, 8\}$}\label{fig:8config}
\end{center}
\end{figure}

   Note that in each of the three cases there exists an orbit of triangles
satisfying the conditions to apply the collapsing operation, namely the orbit
containing triangles with flags labelled $1$, $2$, $3$, $10$, $11$, $14$ for
cases (B) and (C), and labelled $2$, $3$, $9$, $8$, $7$, $6$ for
case (D).
By applying the collapsing operation we effectively delete all flags on the deleted
triangles preserving the remaining flags as well as
all symmetries of the map. It may be the case that we gain some extra symmetries and
therefore we obtain a $k$-orbit map for some divisor $k$ of
$8$. The new map will be vertex transitive and
admissible to the respective types defined by the local
arrangements of flags in Figure \ref{fig:8config}. Observe that cases (C)
and (D) are equal up to relabeling.

\begin{figure}
\begin{center}
\includegraphics[width=6.5cm, height=3cm]{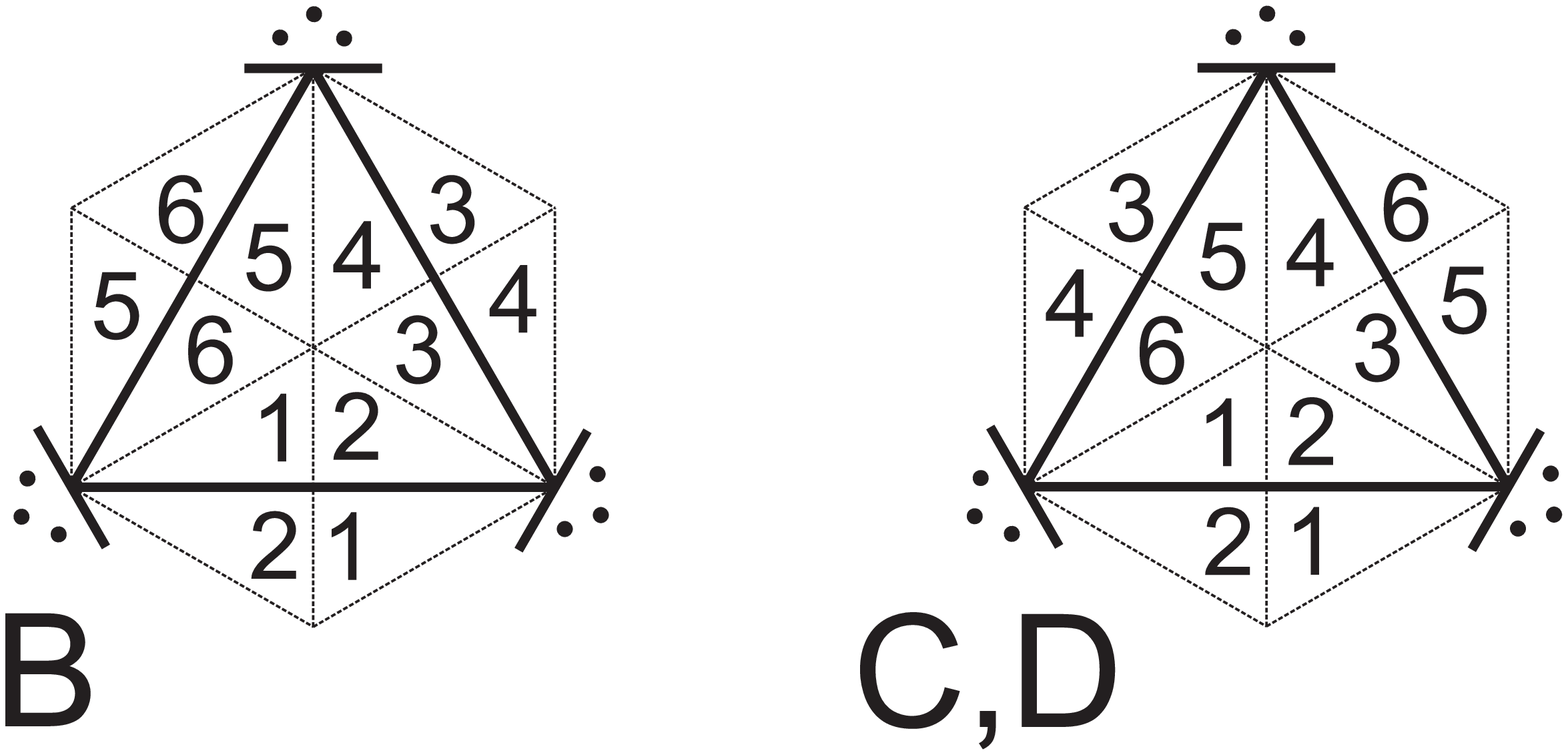}
\caption{Local configurations of $6$-orbit vertex-transitive maps with type
$\{3, 9\}$}\label{fig:6config}
\end{center}
\end{figure}

   All types in Figure \ref{fig:8config} admit the rebelting operation. The
triangles with labels $12$ and $13$ are deleted while the belts formed by the remaining
triangles are rearranged to obtain a map with Schl\"afli type $\{3, 9\}$.
Therefore the new maps
are $k$-orbit maps for some divisor $k$ of $6$ and are admissible with the local
arrangements of flags in Figure \ref{fig:6config}, where the flags are arbitrarily
labelled in the set $\{1, \dots, 6\}$. Note that in all cases the
new map is not only vertex-transitive, but also face-transitive. Once more
the symmetries of the original $14$-orbit map are preserved.

\begin{figure}
\begin{center}
\includegraphics[width=6.5cm, height=3cm]{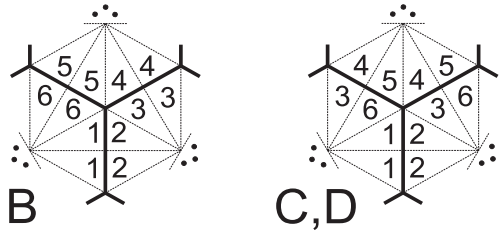}
\caption{Local configurations of $6$-orbit $3$-valent vertex-transitive
maps}\label{fig:6config2}
\end{center}
\end{figure}

   Applying the dual to a map admissible with one of the types in Figure \ref{fig:6config}
we obtain a vertex- and face-transitive map with Schl\"afli type $\{9, 3\}$. Applying the
Petrial operation to the latter yields a $3$-valent map whose Petrie polygons have length
$9$ and must be admissible with the types described by the local configuration of flags
in Figure \ref{fig:6config2}.

\begin{figure}
\begin{center}
\includegraphics[width=5.5cm, height=3cm]{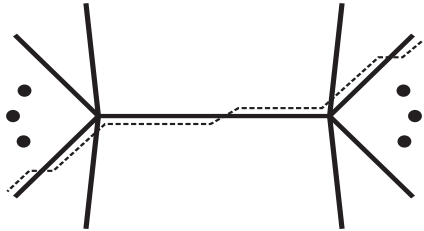}
\caption{$2$-zigzag of a map}\label{fig:2zigzag}
\end{center}
\end{figure}

   Each flag labelled $3$ or $6$ on Figure \ref{fig:6config2} can be glued with its
$1$- and $2$-adjacent flags to assemble new flags with vertices on the centres of faces
with flags labelled $4$ and $5$, midpoints of edges on the midpoints of edges with flags
labelled $1$ and $2$ and centres of faces on the centres of faces with flags labelled
$1$, $2$, $3$ and $6$. This implies that any map $M'$ admissible with any of
the types described by the local configurations in Figure \ref{fig:6config2} is
the truncation of a map $M$ admissible with type $2_{01}$ (case (B)) or with type
$2_0$ (cases (C) and (D)). Note that edges with flags labelled $1$ and $2$ are the
hereditary edges of $M'$ and appear every
three steps in the Petrie polygons. We recall that the edges of $M$ correspond
to the hereditary edges of $M'$. Therefore the Petrie polygons with
length $9$ of $M'$ translate into paths with length $3$ in $M$. These paths are
$2$-zigzags in the sense of the proof of Lemma 7B11 in \cite{ARP}, that is, they traverse
an edge with a given local orientation (left or right) and continue by taking the
second edge in that direction but changing the local orientation
(see Figure \ref{fig:2zigzag}). In particular, $M$ is a map on a non-orientable surface.

   By reversing the process described above we can construct all vertex-transitive maps with
Schl\"afli type $\{3, 7\}$ which are neither regular nor chiral by
determining all maps of type $2_{01}$ satisfying the relation
$(\rho_0 \rho_1 \rho_{212})^3 = id$
and all maps of type $2_0$ satisfying the relation $(\rho_0 \sigma_{12}^2)^3 = id$.
It is easy to verify that these two relations are the ones determining the length $3$
of the $2$-zigzags on maps of the given admissibility. Conversely, the truncation
of any map with $2$-zigzags of size $3$ has Petrie polygons with length $9$.

   The following algorithm produces a
vertex-transitive map $\bar{M}$ with Schl\"afli type $\{3, 7\}$ from a map $M$
which is either $2_{01}$-admissible and satisfies the relation
$(\rho_0 \rho_1 \rho_{212})^3 = id$ or $2_0$-admissible and satisfies the relation
$(\rho_0 \sigma_{12}^2)^3 = id$.

\begin{enumerate}
\item Truncate $M$ to obtain a $3$-valent map $M' := Tr(M)$ with Petrie polygons of length $9$.
   Let $E_1$ be the set of hereditary edges of $M'$
\item Take the petrial $Pe(M')$ of $M'$ to obtain a map with Schl\"afli type $\{9, 3\}$.
   Since the vertex and edge set of $M'$ and $Pe(M')$ coincide we may think of $E_1$
   as a set of edges of $Pe(M')$. It is not hard to verify from the type graph
   that, if $M$ is $2_{01}$-admissible, the resulting surface is orientable.
\item Take the dual of $Pe(M')$ to obtain a map $M'' := Du(Pe(M'))$ with Schl\"afli
   type $\{3, 9\}$. We refer by $E_1'$ the set of edges intersecting those of
   $E_1$. Then the faces of $M''$ can be arranged in belts of triangles determined by the
   edges in $E_1'$.
\item Take the inverse of the rebelting operation on $M''$ with respect to the belts
   described in the previous item. This produces a map $M^{(3)}$ with Schl\"afli type
   $\{3, 8\}$ whose faces can be divided in those arranged on belts
   and those arising from the vertices of $M''$
\item \begin{itemize}
      \item [(a)] If $M$ is $2_{01}$-admissible, to construct a map $\bar{M}$
         admissible with the type represented by the flag arrangement (B) in
         Figure \ref{fig:14config}, apply the inverse of the collapsing operation
         in the following way. Label the flags of $M^{(3)}$ as indicated in Figure
         \ref{fig:8config} (B). This can be done by labelling the flags of
         $M'$ and $M''$ as indicated in Figures \ref{fig:6config2} B and
         \ref{fig:6config} B respectively as intermediate steps.
         For each flag $\Phi$ labelled $k$ for some $k \in \{12, 13\}$ consider its
         defining vertex $v_\Phi$ and edge $e_\Phi$. Let $e_\Phi^{op}$ be the edge
         incident to $v_\Phi$ other than $e_\Phi$ delimiting a flag labelled
         $k$ with a vertex at $v_\Phi$. Split the vertex $v_\Phi$ into
         two vertices $v_1$ and $v_2$ while splitting each of $e_\Phi$ and
         $e_\Phi^{op}$ in two edges, one of them incident to $v_1$ and the
         other to $v_2$. Add an edge between $v_1$ and $v_2$.
      \item [(b)] If $M$ is $2_0$-admissible, to construct a map $\bar{M}$
         admissible with the type represented by the flag arrangement (C) in
         Figure \ref{fig:14config}, apply the inverse of the collapsing operation
         repeating the steps in item (a) replacing the
         labelling of the flags by the one in Figure \ref{fig:8config} (C).
      \item [(c)] Finally, if $M$ is $2_0$-admissible, to construct a map $\bar{M}$
         admissible with the type represented by the flag arrangement (D) in
         Figure \ref{fig:14config}, apply the inverse of the collapsing operation
         in a similar way to that described in item (a), replacing the
         labelling of the flags by the one in Figure \ref{fig:8config} (D), and
         considering all flags labelled $k \in \{1, 10\}$ instead of those with $k \in \{12, 13\}$.
      \end{itemize}
\end{enumerate}

   Let $M$ be a $2_{01}$-admissible (resp. $2_0$-admissible) map with $2_{01}$-admissible
(resp. $2_0$-admissible) group $G$, and $\bar{M}$ be the map obtained from
$M$ by the algorithm above. Since the Petrial and dual operations preserve
the group, and the truncation operation only preserves or increases it,
$G$ is also a subgroup of the automorphism group of the map obtained from $M$ by the
three first steps of the algorithm. Furthermore, the last two steps of
the algorithm effectively add some orbits of flags while preserving the action of
$G$ on the old (and new) flags. Therefore $G$ is a subgroup of the map obtained by
the algorithm above.

   Observe that $M$ and $\bar{M}$ contain $2|G|$ and $14|G|$ flags respectively,
and therefore, $|G|/2$ and $7|G|/2$ edges respectively.
This implies that if $\Gamma(\bar{M}) = G$ then $G$ is a $14$-orbit map, whereas
if $G$ is respectively an index $7$ or $14$ subgroup of $\Gamma(\bar{M})$ then
$\bar{M}$ is chiral or regular. By Proposition \ref{prop:7orbit}, $G$ cannot be an
index $2$ subgroup of $\Gamma(\bar{M})$.

   The Euler characteristic of the surface of
$\bar{M}$ is $v - e + f = 14|G|(1/14 - 1/4 + 1/6) = -|G|/6$ which is one third
of the number of edges of $M$. Consequently, to obtain all vertex-transitive maps
on all surfaces with Euler characteristic $\chi$
with Schl\"afli type $\{3, 7\}$ which are neither regular nor chiral it suffices
to determine all maps with $3\chi$ edges
satisfying the conditions required for the algorithm above.

   Note, however, that a given $2_{01}$- or $2_0$-admissible map yields two vertex-transitive maps
with Schl\"afli type $\{3,7\}$ whenever the two possible choices of $k$ in step 5 of the algorithm are essentially different. If these two choices of $k$ produce isomorphic maps, then there is an extra symmetry of the map obtained after step 4. By tracing back this symmetry to $M$ via the operations, we see that the two choices of $k$ are isomorphic if and only if the map $M$ is regular.

   Furthermore, a given group $G$ generated by three involutions determines three
$2_{01}$-admissible maps depending on the choice of generator to play the role of $\rho_0$. However, in some cases some of these maps may be isomorphic. Consequently, a group $G$ generated by three involutions whose product has order 3 determines $6/t$ maps where $t \in \{1,2,3,6\}$ is the number of isomorphisms of $G$ induced by permutations of its three involutory generators.

   To conclude this section we present a theorem that shows that there are infinitely many surfaces that
do not admit non-regular vertex-transitive maps with Schl\"afli type $\{3, 7\}$.

\begin{theorem}\label{last}
   Let $M$ be a vertex-transitive map with Schl\"afli type $\{3, 7\}$ on a surface with Euler
characteristic $-p$ for $p$ odd prime. Then $M$ must be regular.
\end{theorem}

\begin{proof}
   The discussion in Section \ref{section:3773} shows that $M$ must be either regular, chiral
or $14$-orbit map.

   Surfaces with odd Euler characteristic are non-orientable and therefore $M$ cannot be chiral
or admissible with the flag-configuration in Figure \ref{fig:14config} B.

   It remains to show that $M$ cannot be a $14$-orbit map admissible with the flag-arrangements
in Figure \ref{fig:14config} C, D. Assume to the contrary that such a map exists,
then it can be obtained by the algorithm in
Section \ref{section:3773} from a $2_0$-admissible map with $3p$ edges and
$2_0$-admissible group with $6p$ elements. The case $p=3$ was
discarded by an exhaustive search. Now assume that $p \ge 5$.
According to \cite[Table 16.3]{conway} the only groups with order $6p$ are $\mathbb{Z}_{6p}$,
$\mathbb{Z}_p \times D_3$, $\mathbb{Z}_3 \times D_p$, $D_{3p}$,
$(\mathbb{Z}_p : \mathbb{Z}_q) \times \mathbb{Z}_2$ and $\mathbb{Z}_p : \mathbb{Z}_6$. The last
two occur only if and only if $p \equiv 1 (\mbox{mod 6})$. It is not hard to verify that if a
generating set of any of these groups contains only an involution $a$ and another element $b$ then $ab^2$ cannot have order $3$. This implies that there are no $2_0$-admissible maps with $3p$ edges
and the theorem holds.
\end{proof}

   The regular maps N15.1 and N147.1 in \cite{conderatlas} lie on surfaces with Euler characteristic
$-17$ and $-149$ respectively, both negatives of prime numbers showing that Theorem \ref{last} is false if we
omit the assumption of non-regularity of $M$.

\section{Examples}\label{sec:exam}

   In this section we present in detail three vertex-transitive maps with Schl\"afli type
$\{3, 7\}$ obtained by the algorithm in Section \ref{section:3773}, one for each
of the items $5$(a), $5$(b) and $5$(c). We also determine the number of vertex-transitive maps with Schl\"afli type $\{3,7\}$ on surfaces with Euler characteristic $-2$.

\begin{figure}
\begin{center}
\includegraphics[width=11.7cm, height=15cm]{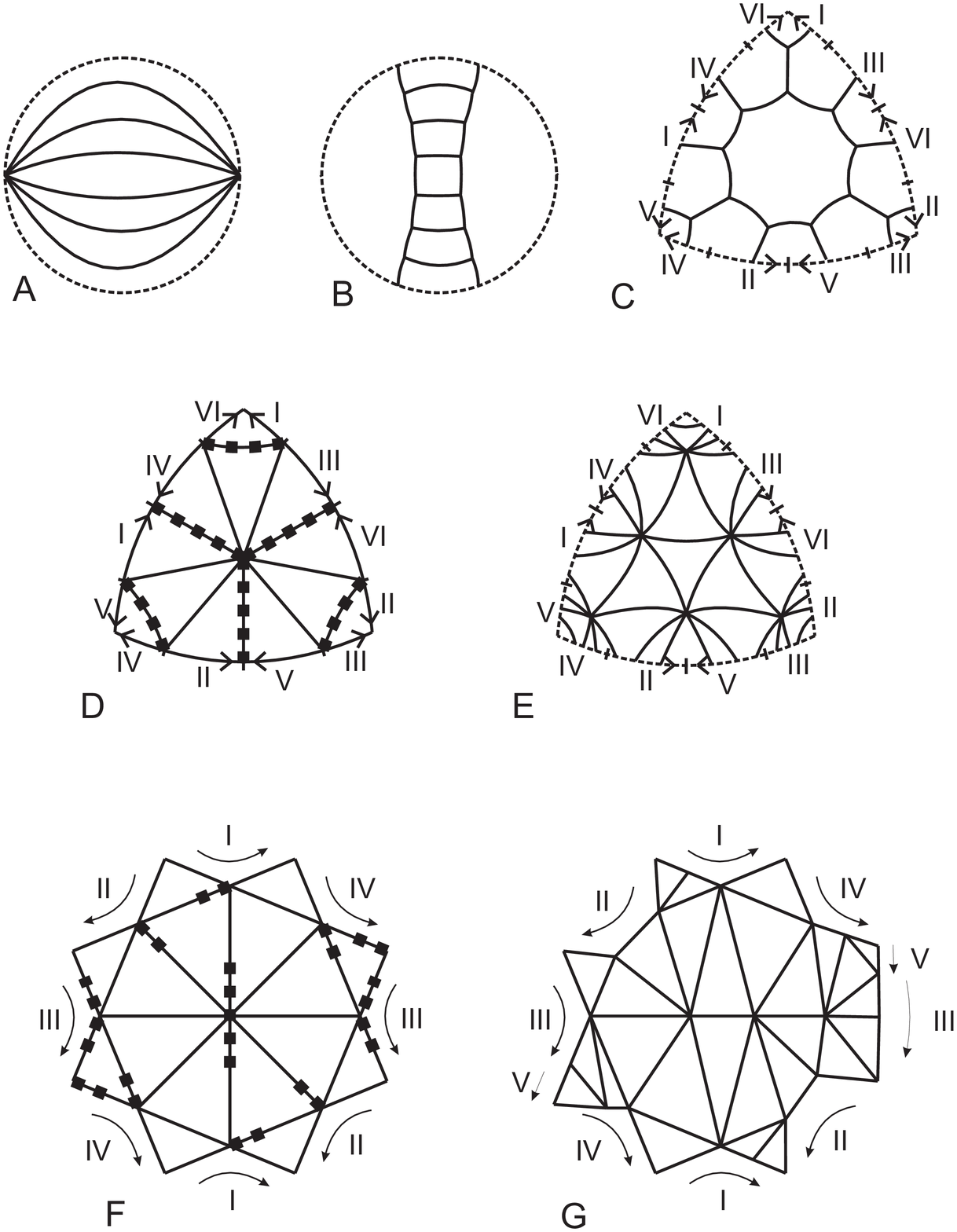}
\caption{Construction of a vertex-transitive map with Schl\"afli type $\{3, 7\}$
on an orientable surface with Euler characteristic $-2$}\label{fig:orient}
\end{center}
\end{figure}

   Figure \ref{fig:orient} illustrates the steps to obtain a vertex-transitive map
with Schl\"afli type $\{3, 7\}$ admissible with the local arrangements of flags in Figure
\ref{fig:14config} B from the $2_{01}$-admissible map in the projective plane
consisting of six $2$-gons (Figure \ref{fig:orient} A). Figure \ref{fig:orient} B
shows a hemi-prism which is the truncation of the map in Figure
\ref{fig:orient} A. The Petrial of the hemi-prism
is illustrated in Figure \ref{fig:orient} C. The surface now is orientable and has Euler
characteristic $-2$ (genus $2$). Figure \ref{fig:orient} D shows the dual of the map in
Figure \ref{fig:orient}. The edges in $E_1$ are indicated in thick dotted lines (see the
algorithm in Section \ref{section:3773}). The inverse of the rebelting
operation applied to the map in Figure \ref{fig:orient} D is shown in Figure \ref{fig:orient} E.
This is the regular map $\{3, 8\}*96$ in \cite{sma} or R2.1 in \cite{conderatlas}
(see also \cite[Figure 8]{alenasiadan}).
Figure \ref{fig:orient} F shows another presentation of the map in Figure \ref{fig:orient}.
The edges in thick dotted lines indicate the places where the inverse of the collapsing
operation must be applied. Finally, Figure \ref{fig:orient} G shows the vertex-transitive map
we were looking for.

\begin{figure}
\begin{center}
\includegraphics[width=11.7cm, height=5.5cm]{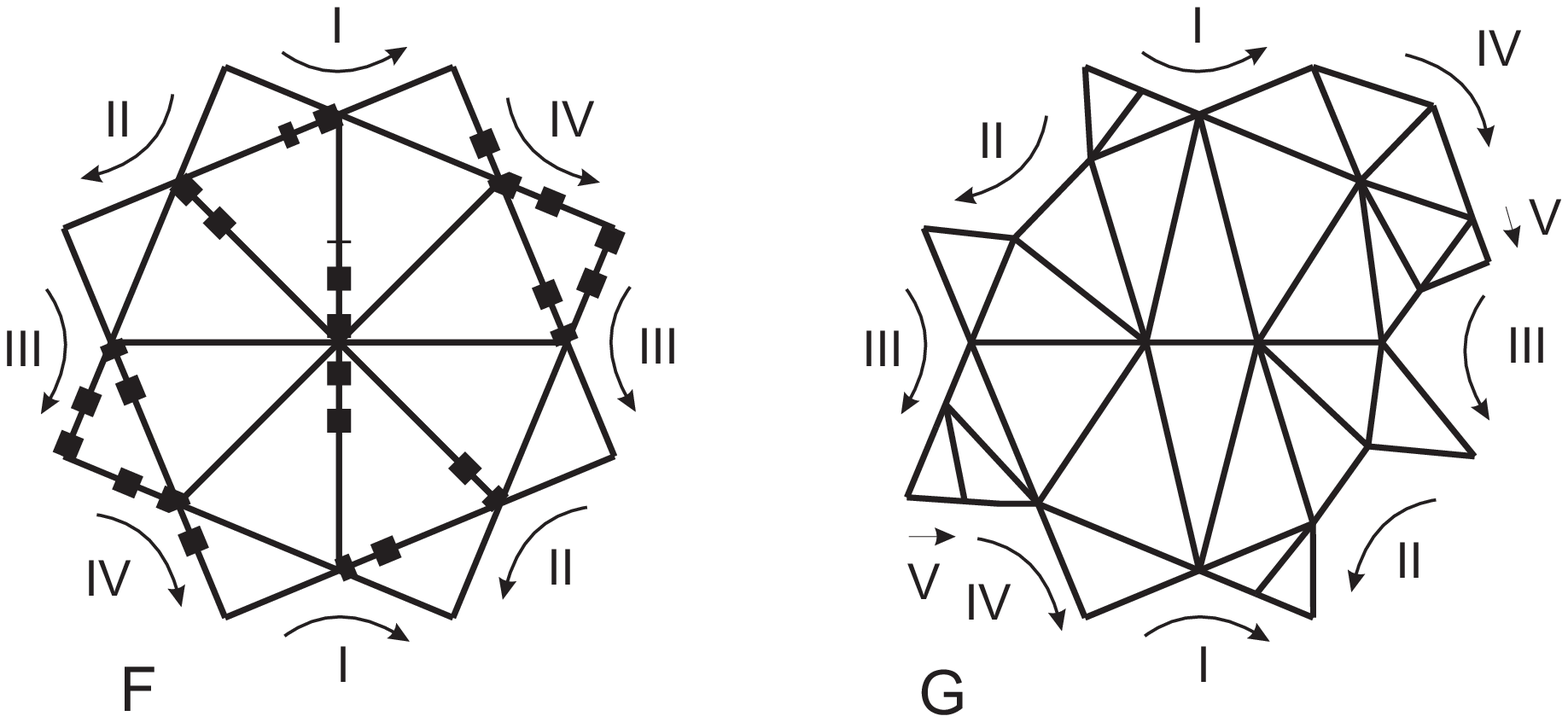}
\caption{Construction of a vertex-transitive map with Schl\"afli type $\{3, 7\}$
on an orientable surface with Euler characteristic $-2$}\label{fig:nonorient1}
\end{center}
\end{figure}

   Figure \ref{fig:nonorient1} illustrates the last steps to obtain a vertex-transitive map
with Schl\"afli type $\{3, 7\}$ admissible with the local arrangements of flags in Figure
\ref{fig:14config} C from the $2_0$-admissible map in the projective plane
with six $2$-gons (Figure \ref{fig:orient} A). Since steps 1 to 4 of the algorithm are the same as for $2_{01}$-admissible maps, these steps are illustrated also by Figure \ref{fig:orient} B, C, D, E. The edges in dark dotted lines in Figure \ref{fig:nonorient1} F indicate where the
inverse of the collapsing operation should be applied. Note that these edges are not
equivalent to the edges in \ref{fig:orient} F. Figure \ref{fig:nonorient1} G corresponds
to the desired map with Schl\"afli type $\{3, 7\}$.

\begin{figure}
\begin{center}
\includegraphics[width=10cm, height=6.7cm]{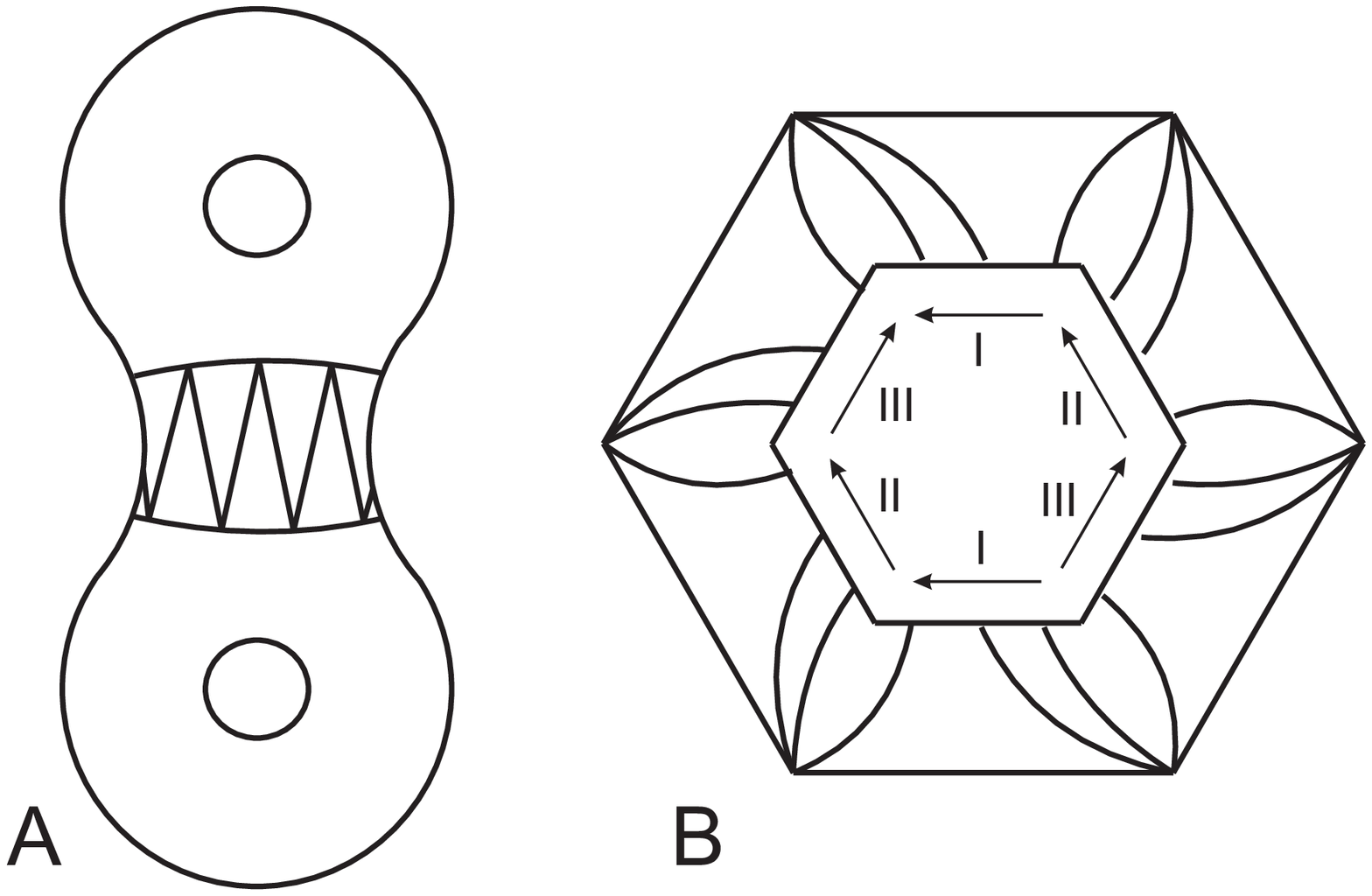}
\caption{Vertex-transitive map with Schl\"afli type $\{3, 7\}$
on an orientable surface with Euler characteristic $-2$}\label{fig:rotvsrefl}
\end{center}
\end{figure}

   The two maps obtained so far can be visualized as follows. Consider a belt formed by $12$
triangles and embed it on an orientable surface of Euler characteristic $-2$ (genus $2$) in such a way that its removal
divides the surface in two tori, each of them missing a disk (see Figure \ref{fig:rotvsrefl}).
On each of these tori embed the map in Figure \ref{fig:rotvsrefl} B, where the outer
hexagon is identified with the border of the belt and the inner hexagon is identified as
indicated by the arrows to form the torus. If both tori are attached to the belt following
the same orientation we obtain the map in Figure \ref{fig:orient} G, whereas if the tori are
attached with opposite orientations we obtain the map in Figure \ref{fig:nonorient1} G.

\begin{figure}
\begin{center}
\includegraphics[width=9cm, height=4.9cm]{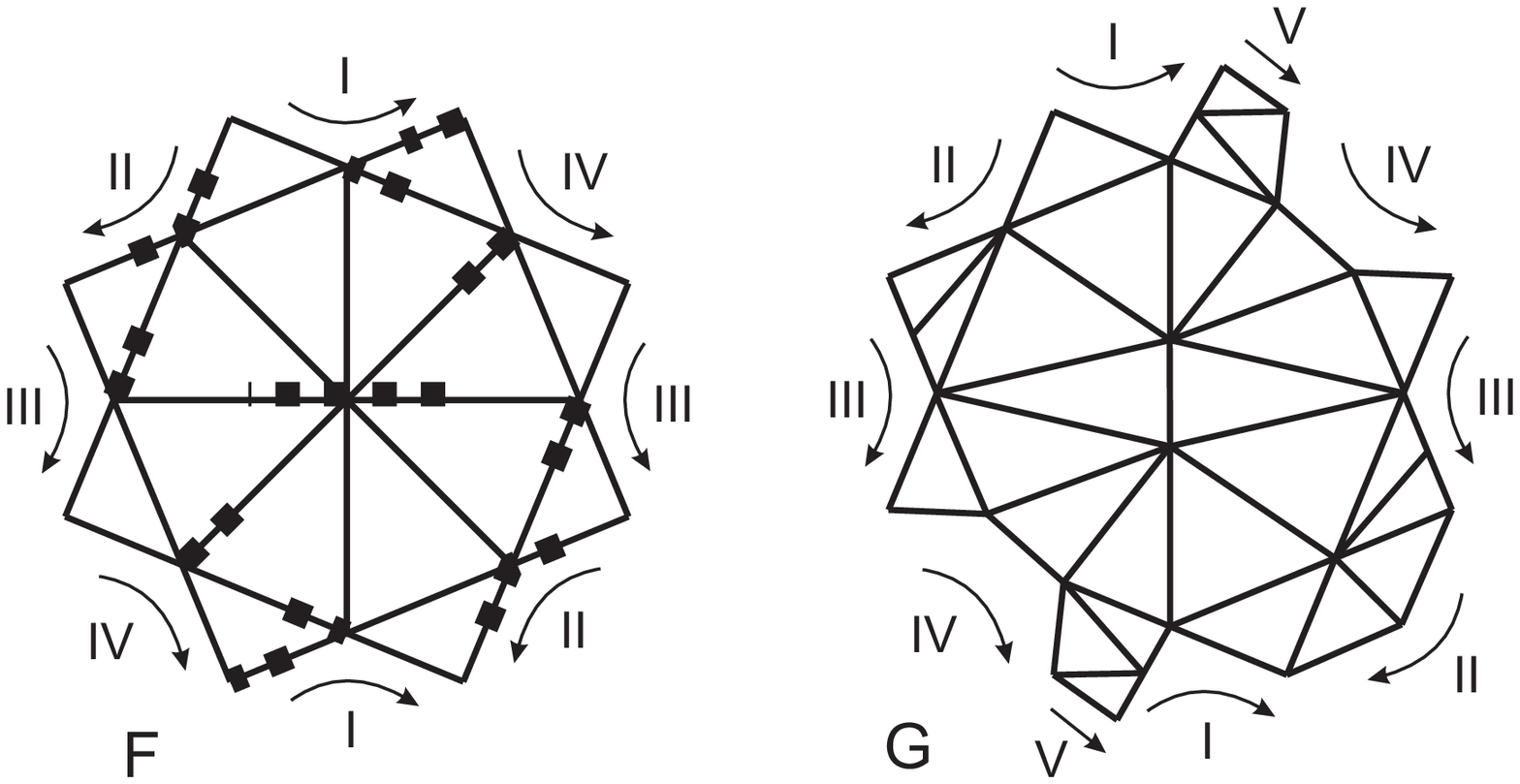}
\caption{Construction of a vertex-transitive map with Schl\"afli type $\{3, 7\}$
on an orientable surface with Euler characteristic $-2$}\label{fig:nonorient2}
\end{center}
\end{figure}

   To obtain a vertex-transitive map
with Schl\"afli type $\{3, 7\}$ admissible with the local arrangements of flags in Figure
\ref{fig:14config} D from the $2_0$-admissible map in the projective plane
with six $2$-gons in Figure \ref{fig:orient} A we follow the steps corresponding to
the pictures in Figure \ref{fig:orient} A--E. Figure \ref{fig:nonorient2} F shows the
edges where the inverse of the collapsing operation should be applied and the desired
map is the one in Figure \ref{fig:nonorient2} G.

   The map in Figure \ref{fig:orient} A and the $2_{01}$-compatible map on the projective plane
consisting of a hexagon and three $2$-gons are the only $2_{01}$-compatible maps with 6 edges. Since the map in Figure \ref{fig:orient} A is regular, it only yields one vertex-transitive map with Schl\"afli type $\{3,7\}$. The remaining $2_{01}$-compatible map with 6 edges is not regular and hence it yields two non-isomorphic vertex-transitive maps with Schl\"afli type $\{3,7\}$ on an orientable surface with Euler characteristic $-2$. Equivalently, the dihedral group with 12 elements $D_6$ is generated by three involutions whose product has order $3$; since there is precisely one non-trivial automorphism of $D_6$ permuting these involutory generators, there are precisely three vertex-transitive maps with Schl\"afli type $\{3,7\}$ associated with this group. These three maps correspond to those arising from the map in Figure \ref{fig:orient} A and the map on the projective plane with a hexagon and three $2$-gons.

   The map in Figure \ref{fig:orient} A and the hemicube are the only $2_0$-admissible maps with 6
edges. Since both of them are regular, the maps in Figures \ref{fig:nonorient1} G and \ref{fig:nonorient2} G, and the two maps arising from the hemicube are the only vertex-transitive maps with Schl\"afli type $\{3,7\}$ arising from $2_0$-compatible maps with 6 edges. The maps in Figures \ref{fig:nonorient1} G and \ref{fig:nonorient2} G lie on an orientable surface, whereas the maps arising from the hemicube lie on a non-orientable surface. In both cases the Euler characteristic is $-2$. This implies that there are precisely 5 vertex-transitive maps with Schl\"afli type $\{3,7\}$ on orientable surfaces with Euler characteristic $-2$, and precisely 2 vertex-transitive maps with Schl\"afli type $\{3,7\}$ on non-orientable surfaces with Euler characteristic $-2$.

\section{Acknowledgements}

   The author wishes to acknowledge Mark Mixer for helpful discussion, and the anonymous referees for useful suggestions.

\appendix
\section{Atlas of vertex-transitive maps with Schl\"afli type $\{3,7\}$ by genus}

   Here we present information about all vertex-transitive maps with Schl\"afli type $\{3, 7\}$ on surfaces
with Euler characteristic $-1 \ge \chi \ge -40$.

   Each vertex-transitive map with Schl\"afli type $\{3, 7\}$ on a surface with Euler
characteristic $\chi$ which is neither regular nor chiral can be obtained from either a
$2_{01}$-admissible map with $-3\chi$ edges satisfying the relation
$(\rho_0 \rho_1 \rho_{212})^3 = id$, or from a $2_0$-admissible map with $-3\chi$ edges
satisfying the relation $(\rho_0 \sigma_{12})^3 = id$. This is equivalent to find all
groups with order $-6\chi$ generated by $3$ involutions whose product has order $3$ and
all groups with order $-6\chi$ generated by an involution and a non-involution such that
the product of the former with the square of the latter has order $3$.

   In the lists of regular and chiral maps \cite{conderatlas} there is $1$ chiral
map and $8$ regular maps with Schl\"afli type $\{3, 7\}$, five of the regular
maps lie on orientable surfaces and the remaining three on non-orientable surfaces.
The regular map $R3.1$ and the chiral map
$C17.1$ can be obtained from the algorithm with the step $5$(a), the regular
map $R14.1$ can be obtained from the algorithm with the step $5$(b), and the
regular map $R14.3$ can be obtained from the algorithm with the step $5$(c).
The remaining maps in these lists cannot be obtained by the construction above.

   As explained at the end of Section \ref{section:3773}, a given $2_{01}$-admissible group $G$ together with its generators
may determine one, two, three or six distinct vertex-transitive maps with Schl\"afli type $\{3,7\}$, depending on the number of permutations of the involutory generators extending to automorphisms of $G$. On the other hand, a given $2_0$-admissible group $G$ with given generators $\rho_0$ and $\sigma_{12}$ may determine one or two vertex-transitive maps with Schl\"afli type $\{3,7\}$ depending on whether there is an automorphism of $G$ fixing $\rho_0$ and mapping $\sigma_{12}$ to $\sigma_{12}^{-1}$ (that is, the corresponging $2_{01}$-admissible map is regular). All these possibilities occur in the atlas with the exception of a $2_{01}$-admissible group yielding only two vertex-transitive maps with Schl\"afli type $\{3,7\}$.

   The following list contains all pairs consisting of a group and a generating set which are subgroups of the automorphism group of a vertex-transitive maps with Schl\"afli type $\{3,7\}$ on surfaces with Euler characteristic $-1 \ge \chi \ge -40$. With each such group we indicate one of the corresponding $2_0$- or $2_{01}$-admissible maps according to the lists in \cite{conderatlas}. Sometimes the $2_{01}$-admissible map is not regular but can be obtained from a regular map by doubling every edge and interpreting the induced $2$-gons as faces. We also mention how many vertex-transitive maps with Schl\"afli type $\{3,7\}$ are associated to each group. Notice that for several genera $\chi$ there are no vertex-transitive maps with Schl\"afli type $\{3,7\}$ on surfaces with Euler characteristic $\chi$.

{\bf Automorphism groups $G$ of vertex-transitive maps with Schl\"afli type $\{3,7\}$ arising from $2_{01}$-admissible maps.}

$\blacktriangleright \,\,$ $\chi = -2$ ($|G|=12$)

\begin{itemize}
   \item $\rho_0 = (1,2)(3,4)(5,6)$,\\
$\rho_1 = (2,3)(4,5)$,\\
$\rho_{212} = (1,6)(2,5)(3,4)$.\\
{\em $2_{01}$-admissible map}: regular map with six $2$-gons on the projective plane.\\
This group yields $3$ different $14$-orbit vertex-transitive maps with Schl\"afli type $\{3,7\}$.
\end{itemize}

$\blacktriangleright \,\,$ $\chi = -4$ ($|G|=24$)

\begin{itemize}
   \item $\rho_0 = (1,2)$,\\
$\rho_1 = (2,3)$,\\
$\rho_{212} = (1,2)(3,4)$.\\
{\em $2_{01}$-admissible map}: hemicuboctahedron on the projective plane.\\
This group yields $6$ different $14$-orbit vertex-transitive maps with Schl\"afli type $\{3,7\}$.

   \item $\rho_0 = (1,2)$,\\
$\rho_1 = (2,3)$,\\
$\rho_{212} = (1,3)(2,4)$.\\
{\em $2_{01}$-admissible map}: N4.1 in \cite{conderatlas}.\\
This group yields $2$ different $14$-orbit vertex-transitive maps and $1$ regular map with Schl\"afli type $\{3,7\}$.
\end{itemize}

$\blacktriangleright \,\,$ $\chi = -6$ ($|G|=36$)

\begin{itemize}
   \item $\rho_0 = (1, 2)(3, 4)(5, 6)$,\\
$\rho_1 = (2, 3)(4, 5)$,\\
$\rho_{212} = (1, 3)(2, 4)(5, 6)$.\\% 
{\em $2_{01}$-admissible map}: dual of N5.2 in \cite{conderatlas}.\\
This group yields $3$ different $14$-orbit vertex-transitive maps with Schl\"afli type $\{3,7\}$.
\end{itemize}

$\blacktriangleright \,\,$ $\chi = -8$ ($|G|=48$)

\begin{itemize}
   \item $\rho_0 = (2, 4)(6, 7)$,\\
   $\rho_1 =  (1, 2)(3, 6)(4, 5)(7, 8)$,\\%
   $\rho_{212} = (1, 3)(2, 5)(4, 7)(6, 8)$.\\
{\em $2_{01}$-admissible map}: N10.2 in \cite{conderatlas}.\\
This group yields $3$ different $14$-orbit vertex-transitive maps with Schl\"afli type $\{3,7\}$.

   \item $\rho_0 = (1, 2)(3, 5)(4, 6)(7, 8)$,\\
   $\rho_1 = (2, 4)(5, 7)$,\\
   $\rho_{212} = (1, 3)(2, 5)(4, 6)(7, 8)$.\\% N4.2 and its dual with double edges, plus the medial of N4.2
{\em $2_{01}$-admissible map}: N4.2 in \cite{conderatlas} with double edges.\\
This group yields $6$ different $14$-orbit vertex-transitive maps with Schl\"afli type $\{3,7\}$.
\end{itemize}

$\blacktriangleright \,\,$ $\chi = -10$ ($|G|=60$)

\begin{itemize}
   \item $\rho_0 = (3, 4)(5, 6)$,\\
   $\rho_1 = (1, 2)(3, 4)$,\\
   $\rho_{212} = (2, 3)(4, 5)$.\\% 
{\em $2_{01}$-admissible map}: dual of N5.1 in \cite{conderatlas}.\\
This group yields $3$ different $14$-orbit vertex-transitive maps with Schl\"afli type $\{3,7\}$.

   \item $\rho_0 = (1, 2)(5, 6)$,\\
   $\rho_1 = (2, 3)(4, 6)$,\\
   $\rho_{212} = (2, 4)(3, 5)$.\\% 
{\em $2_{01}$-admissible map}: dual of N6.1 in \cite{conderatlas}.\\
This group yields $3$ different $14$-orbit vertex-transitive maps with Schl\"afli type $\{3,7\}$.

   \item $\rho_0 = (2, 3)(5, 6)$,\\
   $\rho_1 =  (1, 2)(4, 5)$,\\
   $\rho_{212} = (2, 4)(3, 5)$.\\% 
{\em $2_{01}$-admissible map}: N10.6 in \cite{conderatlas}.\\
This group yields $3$ different $14$-orbit vertex-transitive maps with Schl\"afli type $\{3,7\}$.
\end{itemize}

$\blacktriangleright \,\,$ $\chi = -16$ ($|G|=96$)

\begin{itemize}
   \item $\rho_0 = (1, 2)(6, 7)$,\\
   $\rho_1 = (1, 3)(2, 4)(5, 7)(6, 8)$,\\
   $\rho_{212} = (1, 3)(2, 5)(4, 6)(7, 8)$.\\% 
{\em $2_{01}$-admissible map}: N10.1 in \cite{conderatlas}.\\
This group yields $3$ different $14$-orbit vertex-transitive maps with Schl\"afli type $\{3,7\}$.

   \item $\rho_0 = (1, 2)(3, 6)(4, 5)(7, 12)(8, 14)(9, 15)(10, 11)(13, 16)$,\\
   $\rho_1 = (2, 4)(3, 7)(5, 9)(6, 10)(8, 12)(15, 16)$,\\
   $\rho_{212} = (1, 3)(2, 5)(4, 8)(6, 11)(7, 13)(9, 16)(10, 15)(12, 14)$.\\% 
{\em $2_{01}$-admissible map}: N34.6 in \cite{conderatlas}.\\
This group yields $3$ different $14$-orbit vertex-transitive maps with Schl\"afli type $\{3,7\}$.

   \item $\rho_0 = (1, 2)(3, 5)(4, 7)(6, 10)(8, 12)(9, 13)(11, 15)(14, 16)$,\\
   $\rho_1 = (2, 4)(3, 6)(5, 9)(7, 8)(10, 14)(12, 15)$,\\
   $\rho_{212} = (1, 3)(2, 5)(4, 8)(6, 11)(7, 12)(9, 13)(10, 15)(14, 16)$.\\% N16.7 and its dual with double edges, plus the medial of N16.7
{\em $2_{01}$-admissible map}: N16.7 in \cite{conderatlas} with double edges.\\
This group yields $6$ different $14$-orbit vertex-transitive maps with Schl\"afli type $\{3,7\}$.
\end{itemize}

$\blacktriangleright \,\,$ $\chi = -18$ ($|G|=108$)

\begin{itemize}
   \item $\rho_0 = (3, 4)(5, 6)(7, 9)$,\\
   $\rho_1 = (2, 3)(4, 5)(6, 8)$,\\
   $\rho_{212} = (1, 2)(3, 4)(5, 7)(6, 9)$.\\% 
{\em $2_{01}$-admissible map}: dual of N11.1 in \cite{conderatlas}.\\
This group yields $3$ different $14$-orbit vertex-transitive maps with Schl\"afli type $\{3,7\}$.

   \item $\rho_0 = (1, 2)(3, 6)(4, 8)(5, 9)(7, 12)(10, 11)$,\\
   $\rho_1 = (1, 3)(2, 4)(5, 10)(6, 11)(7, 9)(8, 12)$,\\
   $\rho_{212} = (1, 4)(2, 5)(3, 7)(6, 10)(8, 11)(9, 12)$.\\% 
{\em $2_{01}$-admissible map}: N29.3 in \cite{conderatlas}.\\
This group yields $1$ different $14$-orbit vertex-transitive maps with Schl\"afli type $\{3,7\}$.

   \item $\rho_0 = (2, 3)(4, 7)(6, 8)$,\\
   $\rho_1 = (2, 4)(3, 5)(7, 9)$,\\
   $\rho_{212} = (1, 2)(3, 6)(5, 7)(8, 9)$.\\% N29.4 
{\em $2_{01}$-admissible map}: N29.3 in \cite{conderatlas}.\\
This group yields $3$ different $14$-orbit vertex-transitive maps with Schl\"afli type $\{3,7\}$.

   \item $\rho_0 = (2, 3)(5, 7)(6, 8)$,\\
   $\rho_1 = (2, 4)(3, 5)(7, 9)$,\\
   $\rho_{212} = (1, 2)(3, 6)(4, 5)(8, 9)$.\\% N29.5
{\em $2_{01}$-admissible map}: N29.3 in \cite{conderatlas}.\\
This group yields $3$ different $14$-orbit vertex-transitive maps with Schl\"afli type $\{3,7\}$.
\end{itemize}

$\blacktriangleright \,\,$ $\chi = -20$ ($|G|=120$)

\begin{itemize}
   \item $\rho_0 = (4, 6)(7, 8)(9, 11)(10, 12)$,\\
   $\rho_1 = (1, 2)(3, 5)(4, 7)(6, 8)$,\\
   $\rho_{212} = (1, 3)(2, 4)(5, 8)(6, 9)(7, 10)(11, 12)$.\\%
{\em $2_{01}$-admissible map}: dual of N20.1 in \cite{conderatlas}.\\
This group yields $3$ different $14$-orbit vertex-transitive maps with Schl\"afli type $\{3,7\}$.

   \item $\rho_0 = (1, 2)(3, 5)(4, 7)(6, 8)(9, 12)(10, 11)$,\\
   $\rho_1 = (2, 3)(4, 8)(5, 9)(7, 11)$,\\
   $\rho_{212} = (2, 4)(3, 6)(5, 10)(7, 9)$.\\%
{\em $2_{01}$-admissible map}: N30.8 in \cite{conderatlas}.\\
This group yields $3$ different $14$-orbit vertex-transitive maps with Schl\"afli type $\{3,7\}$.

   \item $\rho_0 = (2, 4)(6, 9)(7, 11)(8, 12)$,\\
   $\rho_1 = (1, 2)(3, 6)(5, 8)(7, 10)$,\\
   $\rho_{212} = (1, 3)(2, 5)(4, 7)(6, 10)(8, 9)(11, 12)$.\\%
{\em $2_{01}$-admissible map}: dual of N30.11 in \cite{conderatlas}.\\
This group yields $3$ different $14$-orbit vertex-transitive maps with Schl\"afli type $\{3,7\}$.

   \item $\rho_0 = (3, 4)(5, 6)(7, 8)(9, 10)$,\\
   $\rho_1 = (1, 2)(3, 4)(5, 7)(6, 8)(9, 10)(11, 12)$,\\
   $\rho_{212} = (2, 3)(4, 5)(7, 9)(10, 11)$.\\% and its dual with double edges, plus the medial of N14.3
{\em $2_{01}$-admissible map}: N14.3 in \cite{conderatlas} with double edges.\\
This group yields $6$ different $14$-orbit vertex-transitive maps with Schl\"afli type $\{3,7\}$.

   \item $\rho_0 = (1, 2)(3, 6)(4, 7)(5, 9)(8, 12)(10, 11)$,\\
   $\rho_1 = (2, 4)(5, 10)(6, 7)(9, 12)$,\\
   $\rho_{212} = (1, 3)(2, 5)(4, 8)(6, 9)(7, 11)(10, 12)$.\\% NPH26.3
{\em $2_{01}$-admissible map}: obtained from the group of the hypermap NPH26.3 in \cite{conderatlas}.\\
This group yields $6$ different $14$-orbit vertex-transitive maps with Schl\"afli type $\{3,7\}$.

   \item $\rho_0 = (2, 3)(5, 7)(6, 9)(8, 12)$,\\
   $\rho_1 = (1, 2)(3, 5)(4, 6)(7, 10)(8, 11)(9, 12)$,\\
   $\rho_{212} = (2, 4)(3, 6)(5, 8)(7, 11)$.\\% NPH34.9
{\em $2_{01}$-admissible map}: obtained from the group of the hypermap NPH34.9 in \cite{conderatlas}.\\
This group yields $6$ different $14$-orbit vertex-transitive maps with Schl\"afli type $\{3,7\}$.
\end{itemize}

$\blacktriangleright \,\,$ $\chi = -24$ ($|G|=144$)

\begin{itemize}

   \item $\rho_0 = (1, 2)(3, 4)(5, 6)(7, 8)(9, 11)(10, 12)$,\\
   $\rho_1 = (3, 4)(5, 7)(6, 8)$,\\
   $\rho_{212} = (2, 3)(4, 5)(6, 9)(7, 10)(11, 12)$.\\% and its dual with double edges, plus the medial of N20.3
{\em $2_{01}$-admissible map}: N20.3 in \cite{conderatlas} with double edges.\\
This group yields $6$ different $14$-orbit vertex-transitive maps with Schl\"afli type $\{3,7\}$.

   \item $\rho_0 = (1, 2)(3, 5)(4, 7)(6, 8)(9, 12)(10, 11)$,\\
   $\rho_1 = (2, 3)(4, 8)(5, 9)(6, 10)(7, 11)$,\\
   $\rho_{212} = (2, 4)(3, 6)(11, 12)$.\\% NPH32.1
{\em $2_{01}$-admissible map}: obtained from the group of the hypermap NPH32.1 in \cite{conderatlas}.\\
This group yields $6$ different $14$-orbit vertex-transitive maps with Schl\"afli type $\{3,7\}$.

   \item $\rho_0 = (1, 2)(3, 5)(4, 6)(7, 9)(8, 11)(10, 12)$,\\
   $\rho_1 = (2, 3)(4, 5)(6, 8)(7, 10)(9, 11)$,\\
   $\rho_{212} = (2, 4)(5, 7)(8, 12)$.\\% NPH38.1
{\em $2_{01}$-admissible map}: obtained from the group of the hypermap NPH84.1 in \cite{conderatlas}.\\
This group yields $6$ different $14$-orbit vertex-transitive maps with Schl\"afli type $\{3,7\}$.
\end{itemize}

$\blacktriangleright \,\,$ $\chi = -32$ ($|G|=192$)

\begin{itemize}

   \item $\rho_0 = (1, 2)(3, 4)(5, 8)(6, 10)(7, 12)(9, 14)(11, 13)(15, 16)$,\\
   $\rho_1 = (1, 3)(2, 5)(4, 7)(6, 11)(8, 12)(9, 15)(10, 14)(13, 16)$,\\
   $\rho_{212} = (1, 4)(2, 6)(3, 7)(5, 9)(8, 13)(10, 12)(11, 16)(14, 15)$.\\% N34.1
{\em $2_{01}$-admissible map}: N34.1 in \cite{conderatlas}.\\
This group yields $3$ different $14$-orbit vertex-transitive maps with Schl\"afli type $\{3,7\}$.

   \item $\rho_0 = (1, 2)(3, 4)(5, 6)(7, 8)$,\\
   $\rho_1 = (2, 3)(4, 6)(5, 7)$,\\
   $\rho_{212} = (2, 4)(3, 5)(6, 7)$.\\% dual of N42.2
{\em $2_{01}$-admissible map}: dual of N42.2 in \cite{conderatlas}.\\
This group yields $3$ different $14$-orbit vertex-transitive maps with Schl\"afli type $\{3,7\}$.

   \item $\rho_0 = (1, 2)(3, 6)(4, 5)(7, 8)$,\\
   $\rho_1 = (2, 4)(5, 7)(6, 8)$,\\
   $\rho_{212} = (1, 3)(2, 5)(6, 7)$.\\% N58.10
{\em $2_{01}$-admissible map}: N58.10 in \cite{conderatlas}.\\
This group yields $3$ different $14$-orbit vertex-transitive maps with Schl\"afli type $\{3,7\}$.

   \item $\rho_0 = (1, 2)(3, 4)(5, 6)(7, 9)(8, 11)(10, 13)(12, 14)(15, 16)$,\\
   $\rho_1 = (2, 3)(4, 5)(6, 8)(7, 10)(12, 15)(13, 14)$,\\
   $\rho_{212} = (3, 4)(5, 7)(6, 9)(8, 12)(11, 14)(15, 16)$.\\% and its dual with double edges, plus the medial of N22.4
{\em $2_{01}$-admissible map}: N22.4 in \cite{conderatlas} with double edges.\\
This group yields $6$ different $14$-orbit vertex-transitive maps with Schl\"afli type $\{3,7\}$.

   \item $\rho_0 = (1, 2)(3, 5)(4, 7)(6, 8)$,\\
   $\rho_1 = (2, 3)(4, 8)(6, 7)$,\\
   $\rho_{212} = (2, 4)(3, 6)(7, 8)$.\\% NPH30.4
{\em $2_{01}$-admissible map}: obtained from the group of the hypermap NPH30.4 in \cite{conderatlas}.\\
This group yields $6$ different $14$-orbit vertex-transitive maps with Schl\"afli type $\{3,7\}$.

   \item $\rho_0 = (1, 2)(3, 4)(5, 7)(6, 8)$,\\
   $\rho_1 = (2, 4)(5, 7)(6, 8)$,\\
   $\rho_{212} = (1, 3)(2, 5)(4, 6)$.\\% NPH46.5
{\em $2_{01}$-admissible map}: obtained from the group of the hypermap NPH46.5 in \cite{conderatlas}.\\
This group yields $6$ different $14$-orbit vertex-transitive maps with Schl\"afli type $\{3,7\}$.

   \item $\rho_0 = (1, 2)(3, 5)(4, 7)(6, 8)(9, 10)(11, 15)(12, 13)(14, 16)$,\\
   $\rho_1 = (2, 3)(4, 8)(5, 9)(7, 11)(10, 14)(13, 15)$,\\
   $\rho_{212} = (2, 4)(3, 6)(5, 10)(8, 12)(9, 13)(11, 16)$.\\% NPH46.6
{\em $2_{01}$-admissible map}: obtained from the group of the hypermap NPH46.6 in \cite{conderatlas}.\\
This group yields $6$ different $14$-orbit vertex-transitive maps with Schl\"afli type $\{3,7\}$.

   \item $\rho_0 = (1, 2)(3, 5)(4, 6)(7, 10)(8, 12)(9, 14)(11, 16)(13, 15)$,\\
   $\rho_1 = (2, 3)(5, 7)(6, 9)(8, 13)(10, 11)(12, 14)$,\\
   $\rho_{212} = (2, 4)(3, 6)(5, 8)(7, 11)(9, 15)(14, 16)$.\\% NPH58.10
{\em $2_{01}$-admissible map}: obtained from the group of the hypermap NPH58.10 in \cite{conderatlas}.\\
This group yields a chiral map map (C17.1 in \cite{conderatlas}) and $5$ different $14$-orbit vertex-transitive maps with Schl\"afli type $\{3,7\}$.
\end{itemize}

$\blacktriangleright \,\,$ $\chi = -40$ ($|G|=240$)

\begin{itemize}

   \item $\rho_0 = (4, 6)(7, 8)(9, 11)(10, 13)(12, 15)(14, 17)(16, 19)(18, 20)$,\\
   $\rho_1 = (1, 2)(3, 5)(4, 7)(6, 8)(9, 12)(10, 14)(11, 15)(13, 17)(16, 20)(18, 19)(21, 23)(22, 24)$,\\
   $\rho_{212} = (1, 3)(2, 4)(5, 8)(6, 9)(7, 10)(11, 13)(12, 16)(14, 18)(15, 17)(19, 21)(20, 22)(23, 24)$.\\% dual of N38.2
{\em $2_{01}$-admissible map}: N38.2 in \cite{conderatlas}.\\
This group yields $3$ different $14$-orbit vertex-transitive maps with Schl\"afli type $\{3,7\}$.

   \item $\rho_0 = (1, 2)(3, 6)(4, 7)(5, 9)(8, 15)(10, 17)(11, 18)(12, 19)(13, 20)(14, 21)(16, 22)(23, 24)$,\\
   $\rho_1 = (2, 4)(5, 10)(6, 11)(7, 13)(9, 16)(12, 15)(18, 23)(19, 21)$,\\
   $\rho_{212} = (1, 3)(2, 5)(4, 8)(6, 12)(7, 14)(9, 13)(10, 15)(11, 17)(16, 21)(18, 22)(19, 23)(20, 24)$.\\% N70.5
{\em $2_{01}$-admissible map}: N70.5 in \cite{conderatlas}.\\
This group yields $3$ different $14$-orbit vertex-transitive maps with Schl\"afli type $\{3,7\}$.

   \item $\rho_0 = (2, 4)(6, 10)(7, 12)(8, 14)(9, 15)(13, 19)(16, 22)(21, 23)$,\\
   $\rho_1 = (1, 2)(3, 6)(4, 7)(5, 9)(8, 11)(10, 16)(12, 17)(13, 20)(14, 19)(15, 21)(18, 23)(22, 24)$,\\
   $\rho_{212} = (1, 3)(2, 5)(4, 8)(6, 11)(7, 13)(9, 10)(12, 18)(14, 15)(16, 23)(17, 24)(19, 21)(20, 22)$.\\% N78.4
{\em $2_{01}$-admissible map}: N78.4 in \cite{conderatlas}.\\
This group yields $3$ different $14$-orbit vertex-transitive maps with Schl\"afli type $\{3,7\}$.

\end{itemize}

{\bf Automorphism groups $G$ of vertex-transitive maps on orientable surfaces with Schl\"afli type $\{3,7\}$ arising from $2_0$-admissible maps.}

$\blacktriangleright \,\,$ $\chi = -2$ ($|G|=12$)

\begin{itemize}
  \item $\rho_0 = (1,3)(2,4)$,\\
$\sigma_{12} = (1,2,3,4)(5,6,7)$.\\
{\em $2_{0}$-admissible map}: regular map with six $2$-gons on the projective plane.\\
This group yields $2$ different $14$-orbit vertex-transitive maps with Schl\"afli type $\{3,7\}$.
\end{itemize}

$\blacktriangleright \,\,$ $\chi = -4$ ($|G|=24$)

\begin{itemize}
  \item $\rho_0 = (1,2)(3,4)$,\\
$\sigma_{12} = (1,2,3)(5,6)$.\\
{\em $2_{0}$-admissible map}: N4.1 in \cite{conderatlas}.\\
This group yields $2$ different $14$-orbit vertex-transitive maps with Schl\"afli type $\{3,7\}$.
\end{itemize}

$\blacktriangleright \,\,$ $\chi = -6$ ($|G|=36$)

\begin{itemize}
  \item $\rho_0 = (2, 3)(4, 5)$,\\
$\sigma_{12} = (1, 2)(3, 4, 6, 5)$.\\ % dual of N5.2
{\em $2_{0}$-admissible map}: dual of N5.2 in \cite{conderatlas}.\\
This group yields $2$ different $14$-orbit vertex-transitive maps with Schl\"afli type $\{3,7\}$.
\end{itemize}

$\blacktriangleright \,\,$ $\chi = -8$ ($|G|=48$)

\begin{itemize}
  \item $\rho_0 = (1, 2)(3, 5)(4, 7)(6, 10)(8, 12)(9, 14)(11, 13)(15, 16)$,\\
$\sigma_{12} = (1, 3, 6, 11)(2, 4, 8, 13, 12, 16, 10, 15, 14, 5, 9, 7)$.\\ % N10.2
{\em $2_{0}$-admissible map}: N10.2 in \cite{conderatlas}.\\
This group yields $2$ different $14$-orbit vertex-transitive maps with Schl\"afli type $\{3,7\}$.
\end{itemize}

$\blacktriangleright \,\,$ $\chi = -16$ ($|G|=96$)

\begin{itemize}
  \item $\rho_0 = (1, 2)(5, 6)$,\\
$\sigma_{12} = (1, 3)(2, 4, 5, 7, 6, 8)$.\\ % N10.1
{\em $2_{0}$-admissible map}: N10.1 in \cite{conderatlas}.\\
This group yields $2$ different $14$-orbit vertex-transitive maps with Schl\"afli type $\{3,7\}$.

  \item $\rho_0 = (1, 2)(3, 5)(4, 7)(6, 10)(8, 12)(9, 14)(11, 13)(15, 19)(16, 21)(17, 20)(18, 24)(22, 27)(23, 28)(25, 29)(26, 30)(31, 32)$,\\
$\sigma_{12} = (1, 3, 6, 11, 17, 23, 29, 32)(2, 4, 8, 13, 12, 18, 25, 24, 14, 5, 9, 15, 20, 19, 26, 31, 30, 21, 10, 16, 22, 28, 27, 7)$.\\ % N34.6
{\em $2_{0}$-admissible map}: N34.6 in \cite{conderatlas}.\\
This group yields $2$ different $14$-orbit vertex-transitive maps with Schl\"afli type $\{3,7\}$.
\end{itemize}

$\blacktriangleright \,\,$ $\chi = -18$ ($|G|=108$)

\begin{itemize}
  \item $\rho_0 = (2, 3)(4, 5)(7, 9)(8, 11)(10, 13)(16, 17)$,\\
$\sigma_{12} = (1, 2)(3, 4, 6, 8)(5, 7, 10, 14)(9, 12)(11, 15, 13, 16)(17, 18)$.\\ % N11.1
{\em $2_{0}$-admissible map}: N11.1 in \cite{conderatlas}.\\
This group yields $2$ different $14$-orbit vertex-transitive maps with Schl\"afli type $\{3,7\}$.

  \item $\rho_0 = (1, 2)(3, 5)(4, 7)(6, 8)(9, 10)(11, 12)$,\\
$\sigma_{12} = (1, 3, 6, 5, 8, 9, 11, 10, 12, 7, 2, 4)$.\\ % N29.3
{\em $2_{0}$-admissible map}: N29.3 in \cite{conderatlas}.\\
This group yields $2$ different $14$-orbit vertex-transitive maps with Schl\"afli type $\{3,7\}$.

  \item $\rho_0 = (1, 2)(4, 6)(5, 7)(8, 11)(12, 15)(13, 17)$,\\
$\sigma_{12} = (1, 3, 5, 8, 7, 10, 13, 15, 17, 18, 2, 4)(6, 9, 12, 16, 11, 14)$.\\ % N29.4
{\em $2_{0}$-admissible map}: N29.4 in \cite{conderatlas}.\\
This group yields $2$ different $14$-orbit vertex-transitive maps with Schl\"afli type $\{3,7\}$.

  \item $\rho_0 = (1, 2)(3, 5)(6, 8)(10, 12)(11, 14)(15, 17)$,\\
$\sigma_{12} = (1, 3, 6, 9, 8, 11, 15, 18, 17, 12, 2, 4)(5, 7, 10, 13, 14, 16)$.\\ % N29.5
{\em $2_{0}$-admissible map}: N29.5 in \cite{conderatlas}.\\
This group yields $2$ different $14$-orbit vertex-transitive maps with Schl\"afli type $\{3,7\}$.
\end{itemize}

$\blacktriangleright \,\,$ $\chi = -20$ ($|G|=120$)

\begin{itemize}
  \item $\rho_0 = (3, 4)(5, 6)$,\\
$\sigma_{12} = (1, 2, 3, 5)$.\\ % Dual of N20.1
{\em $2_{0}$-admissible map}: dual of N20.1 in \cite{conderatlas}.\\
This group yields $2$ different $14$-orbit vertex-transitive maps with Schl\"afli type $\{3,7\}$.

  \item $\rho_0 = (1, 2)(3, 5)(6, 8)(9, 10)$,\\
$\sigma_{12} = (1, 3)(2, 4, 6, 9, 12, 5, 7, 10, 8, 11)$.\\ % N30.8
{\em $2_{0}$-admissible map}: N30.8 in \cite{conderatlas}.\\
This group yields $2$ different $14$-orbit vertex-transitive maps with Schl\"afli type $\{3,7\}$.

  \item $\rho_0 = (1, 2)(3, 4)(5, 7)(6, 9)$,\\
$\sigma_{12} = (2, 3, 5, 8, 7, 9)(4, 6, 10)$.\\ % Dual of N30.11
{\em $2_{0}$-admissible map}: dual of N30.11 in \cite{conderatlas}.\\
This group yields $2$ different $14$-orbit vertex-transitive maps with Schl\"afli type $\{3,7\}$.
\end{itemize}

$\blacktriangleright \,\,$ $\chi = -32$ ($|G|=192$)

\begin{itemize}
  \item $\rho_0 = (1, 2)(3, 4)(5, 7)(6, 9)(8, 11)(10, 12)(13, 14)(15, 16)$,\\
$\sigma_{12} = (1, 3, 5, 8)(2, 4, 6, 10, 9, 12, 7, 11, 13, 15, 14, 16)$.\\ % N34.1
{\em $2_{0}$-admissible map}: N34.1 in \cite{conderatlas}.\\
This group yields $2$ different $14$-orbit vertex-transitive maps with Schl\"afli type $\{3,7\}$.

  \item $\rho_0 = (2, 3)(4, 5)(6, 8)(7, 10)(12, 13)(14, 15)$,\\
$\sigma_{12} = (1, 2)(3, 4, 6, 9, 8, 12)(5, 7, 11, 10, 13, 14)(15, 16)$.\\ % dual of N42.2
{\em $2_{0}$-admissible map}: dual of N42.2 in \cite{conderatlas}.\\
This group yields $2$ different $14$-orbit vertex-transitive maps with Schl\"afli type $\{3,7\}$.

  \item $\rho_0 = (1, 2)(3, 5)(4, 7)(6, 10)(8, 9)(11, 14)(12, 13)(15, 16)$,\\
$\sigma_{12} = (1, 3, 6, 11)(2, 4, 8, 5, 9, 12, 10, 13, 15, 14, 16, 7)$.\\ % N58.10
{\em $2_{0}$-admissible map}: N58.10 in \cite{conderatlas}.\\
This group yields $2$ different $14$-orbit vertex-transitive maps with Schl\"afli type $\{3,7\}$.
\end{itemize}

$\blacktriangleright \,\,$ $\chi = -40$ ($|G|=240$)

\begin{itemize}
  \item $\rho_0 = (3, 4)(5, 6)(7, 8)(11, 12)$,\\
$\sigma_{12} = (1, 2, 3, 5)(6, 7)(8, 9, 10, 11)$.\\ % dual of N38.2
{\em $2_{0}$-admissible map}: dual of N38.2 in \cite{conderatlas}.\\
This group yields $2$ different $14$-orbit vertex-transitive maps with Schl\"afli type $\{3,7\}$.

  \item $\rho_0 = (1, 2)(3, 5)(6, 9)(7, 11)(10, 15)(12, 17)(13, 19)(20, 22)$,\\
$\sigma_{12} = (1, 3, 6, 10)(2, 4, 7, 12, 18, 5, 8, 13, 11, 16, 9, 14, 20, 19, 23, 15, 21, 17, 22, 24)$.\\ % N70.5
{\em $2_{0}$-admissible map}: N70.5 in \cite{conderatlas}.\\
This group yields $2$ different $14$-orbit vertex-transitive maps with Schl\"afli type $\{3,7\}$.

  \item $\rho_0 = (1, 2)(3, 5)(4, 7)(6, 9)(8, 10)(11, 14)(12, 13)(15, 17)(16, 19)(18, 20)$,\\
$\sigma_{12} = (1, 3, 6, 5)(2, 4, 8, 11, 10, 13, 9, 12, 15, 18, 17, 7)(14, 16, 20, 19)$.\\ % N78.4
{\em $2_{0}$-admissible map}: N78.4 in \cite{conderatlas}.\\
This group yields $2$ different $14$-orbit vertex-transitive maps with Schl\"afli type $\{3,7\}$.
\end{itemize}

{\bf Automorphism groups $G$ of vertex-transitive maps with Schl\"afli type $\{3,7\}$ on non-orientable surfaces (arising from $2_0$-admissible maps).}

$\blacktriangleright \,\,$ $\chi = -2$ ($|G|=12$)

\begin{itemize}
  \item $\rho_0 = (1,2)(3,4)$,\\
$\sigma_{12} = (1,2,3)$.\\
{\em $2_{0}$-admissible map}: hemicube (in the projective plane).\\
This group yields $2$ different $14$-orbit vertex-transitive maps with Schl\"afli type $\{3,7\}$.
\end{itemize}

$\blacktriangleright \,\,$ $\chi = -10$ ($|G|=60$)

\begin{itemize}
  \item $\rho_0 = (1, 2)(4, 5)$,\\
$\sigma_{12} = (2, 3, 4, 5, 6)$.\\ % dual of N10.5
{\em $2_{0}$-admissible map}: dual of N10.5 in \cite{conderatlas}.\\
This group yields $2$ different $14$-orbit vertex-transitive maps with Schl\"afli type $\{3,7\}$.
\end{itemize}

$\blacktriangleright \,\,$ $\chi = -26$ ($|G|=156$)

\begin{itemize}
  \item $\rho_0 = (1, 2)(3, 4)(5, 7)(6, 8)(9, 11)(10, 12)$,\\
$\sigma_{12} = (2, 3, 5, 4, 6, 9, 11, 13, 12, 8, 7, 10)$.\\ % dual of C27.7
{\em $2_{0}$-admissible map}: no reference known of this group or of its associated map.\\
This group yields $1$ regular map (R14.1 in \cite{conderatlas}) and $3$ different $14$-orbit vertex-transitive maps with Schl\"afli type $\{3,7\}$.

   \item $\rho_0 = (1, 2)(3, 4)(5, 7)(6, 9)(8, 12)(10, 13)$,\\
$\sigma_{12} = (2, 3, 5, 8, 9, 12, 4, 6, 10, 13, 7, 11)$.\\ % C27.7
{\em $2_{0}$-admissible map}: no reference known of this group or of its associated map.\\
This group yields $1$ regular maps (R14.3 in \cite{conderatlas}) and $3$ different $14$-orbit vertex-transitive maps with Schl\"afli type $\{3,7\}$.
\end{itemize}

$\blacktriangleright \,\,$ $\chi = -28$ ($|G|=168$)

\begin{itemize}
  \item $\rho_0 = (1, 2)(3, 4)(5, 7)(6, 8)$,\\
$\sigma_{12} = (1, 3, 5, 2)(4, 6, 7, 8)$.\\ % dual of N23.2
{\em $2_{0}$-admissible map}: dual of N23.2 in \cite{conderatlas}.\\
This group yields $2$ different $14$-orbit vertex-transitive maps with Schl\"afli type $\{3,7\}$.

  \item $\rho_0 = (1, 2)(3, 4)(5, 7)(6, 8)$,\\
$\sigma_{12} = (2, 3, 5, 4, 6, 7, 8)$.\\ % dual of N41.2
{\em $2_{0}$-admissible map}: dual of N41.2 in \cite{conderatlas}.\\
This group yields $2$ different $14$-orbit vertex-transitive maps with Schl\"afli type $\{3,7\}$.
\end{itemize}

$\blacktriangleright \,\,$ $\chi = -30$ ($|G|=180$)

\begin{itemize}
  \item $\rho_0 = (1, 2)(3, 5)(6, 9)(7, 10)(11, 12)(16, 17)$,\\
$\sigma_{12} = (1, 3, 6)(2, 4, 7, 11, 15, 9, 13, 17, 10, 14, 5, 8, 12, 16, 18)$.\\ % N50.9
{\em $2_{0}$-admissible map}: N50.9 in \cite{conderatlas}.\\
This group yields $2$ different $14$-orbit vertex-transitive maps with Schl\"afli type $\{3,7\}$.
\end{itemize}

   To conclude we show in table \ref{table:table} the number of vertex-transitive maps with Schl\"afli type $\{3,7\}$ which are neither regular nor chiral, on surfaces with Euler characteristic $\chi$ for each $-1 \ge \chi \ge -40$.

\begin{center}
\begin{table}%[h]
\centering
{\small%\tiny
\begin{tabular}{|c|c|c|c|}
\hline
{\textbf{$\chi$}} &
{\textbf{orientable}} &
{\textbf{non-orientable}} &
{\textbf{Total}} \\
\hline
$-2$  & $5$ & 2 & 7  \\
$-4$ & $10$ & 0 & 10 \\
$-6$ & $5$ & 0 & 5 \\
$-8$ & $11$ & 0 & 11 \\
$-10$ & $9$ & 2 & 11 \\
$-16$ & $16$ & 0 & 16 \\
$-18$ & $18$ & 0 & 18 \\
$-20$ & $33$ & 0 & 33 \\
\hline
$-24$ & 18 & 0 & 18\\
$-26$ & 0 & 6 & 6 \\
$-28$ & 0 & 4 & 4 \\
$-30$ & 0 & 2 & 2 \\
$-32$ & 44 & 0 & 44 \\
$-40$ & 15 & 0 & 15 \\
\hline
% \multicolumn{4}{l} { * The labelings of the maps used for the operation belt are described in Section \ref{section8}}
\end{tabular}
}
\caption[]{\label{table:table}Vertex-transitive maps with Schl\"afli type
$\{3, 7\}$ which are neither regular nor chiral on surfaces of Euler characteristic $-1 \ge \chi \ge -40$
}

\end{table}
\end{center}

\end{document}